\RequirePackage{fix-cm}
\documentclass[smallextended]{svjour3}
\smartqed
\usepackage{graphicx,booktabs,amsmath,amssymb}
\usepackage[hidelinks]{hyperref}
\makeatletter
\g@addto@macro{\UrlBreaks}{\do\a\do\b\do\c\do\d\do\e\do\f\do\g\do\h\do\i\do\j\do\k\do\l\do\m\do\n
  \do\o\do\p\do\q\do\r\do\s\do\t\do\u\do\v\do\w\do\x\do\y\do\z\do\A\do\B\do\C\do\D\do\E\do\F
  \do\G\do\H\do\I\do\J\do\K\do\L\do\M\do\N\do\O\do\P\do\Q\do\R\do\S\do\T\do\U\do\V\do\W\do\X
  \do\Y\do\Z\do0\do1\do2\do3\do4\do5\do6\do7\do8\do9}
\makeatother
\usepackage{microtype}
\usepackage{tikz}
\graphicspath{{figs/}}
\journalname{}
\makeatletter\def\makeheadbox{}\makeatother
\date{}

\newcommand{\cl}{C_L}
\newcommand{\cd}{C_D}
\newcommand{\ulp}{\mathrm{ulp}}
\newcommand{\fl}{\mathrm{fl}}
\newcommand{\hi}{\mathrm{hi}}
\newcommand{\lo}{\mathrm{lo}}
\newcommand{\fitwidth}[1]{\resizebox{\ifdim\width>\linewidth\linewidth\else\width\fi}{!}{#1}}

\begin{document}

\title{A wider state, not wider arithmetic, for single-precision time integration}
\titlerunning{A wider state, not wider arithmetic}

\author{Michael F.~Barad}

\institute{M.~F.~Barad \at
  NASA Ames Research Center, Moffett Field, CA 94035, USA \\
  \email{michael.f.barad@nasa.gov} \\
  ORCID: \href{https://orcid.org/0000-0002-7651-0926}{0000-0002-7651-0926}
}

\maketitle

\begin{abstract}
Single precision (FP32) halves the memory and data movement of double precision (FP64), and on
accelerators is often several times faster. Its rounding errors are harmless except in the stored
state. A modified equation for Runge--Kutta methods in Shu--Osher form and linear multistep methods
separates the rounding of the stored state, a forcing proportional to $1/\Delta t$, from truncation
and from right-hand-side rounding, which sets a $\Delta t$-independent floor. We bound this forcing
for strong-stability-preserving schemes and identify two mechanisms that keep its sign fixed, so it
accumulates as a bias: past a turnaround step, refining the step increases the error, in FP64 as in
FP32. For a model oscillator the analysis predicts the floor and the FP32-to-FP64 ratios of
turnaround step and minimum error within $30\%$. Our remedy needs no FP64 arithmetic: each state
variable is stored as the unevaluated sum of two FP32 numbers, updated with error-free
transformations, while every right-hand-side evaluation stays in FP32. Unlike compensated summation,
it needs no rewriting of the scheme. In model problems from geophysical flow, mechanical dynamics
and celestial mechanics, with five Runge--Kutta and multistep schemes, it restores design-order
convergence to the floor. In large-eddy simulations of a high-lift airplane wing on three refined
grids, up to $544$ million cells, with time steps $\Delta t$, $\Delta t/2$ and $\Delta t/4$, FP32
overpredicts the lift by $2.1$, $4.5$ and $7.7\%$, while the two-float lift agrees with FP64 within
$0.3\%$ on every grid, and on graphics processing units the production solver runs $1.6$ to $2.7$
times faster than FP64.
\keywords{Single-precision arithmetic \and Mixed precision \and Round-off error \and Runge--Kutta
methods \and Compensated summation \and Large-eddy simulation}
\subclass{65G50 \and 65L06 \and 65M20 \and 65Y10 \and 76F65}
\end{abstract}

\section{Introduction}\label{sec:intro}

Computational fluid dynamics, like most fields of scientific computing, has traditionally relied on
IEEE double-precision floating-point numbers (FP64) rather than single precision (FP32). That
choice has a cost on any hardware: FP64 doubles the memory occupied by the solution fields and the
data moved per operation, and on processors with vector units it halves the number of values
processed per instruction. On graphics processing units (GPUs), now the dominant platform for
large-scale simulation, the cost is often larger, because their design is increasingly driven by
machine learning, which favors lower precision. On current hardware, FP64 runs at half the FP32
vector rate on the NVIDIA B200 and the AMD MI355X, at a quarter on NVIDIA Rubin, and at one
sixty-fourth on NVIDIA Blackwell Ultra, the AMD MI455X, and workstation GPUs
\cite{nvidia2025blackwell,nvidia2026rubin,amdinstinct}. Fast FP64 is becoming a separate product
for scientific computing, such as the AMD MI430X \cite{amdinstinct}. In the production
compressible-flow solver used in this work, FP32 delivers $1.8$ to $3.1$ times the throughput of
FP64 on the same GPU, and it often determines whether a case fits in GPU memory at all. The analysis
and the remedy developed below, however, do not depend on the hardware: they concern the rounding
of the stored state, which occurs in FP32 on any processor.

The difficulty that this paper addresses can be illustrated with a simple example. At sea level,
air has a density of about $1.225$\,kg\,m$^{-3}$, which FP32 represents with a spacing of about
$1.2\times10^{-7}$\,kg\,m$^{-3}$. A density fluctuation of $10^{-3}$\,kg\,m$^{-3}$, typical of a
flow at Mach $0.2$, spans about eight thousand of these spacings, which is ample to represent the
fluctuation itself. Flow solvers, however, typically store the full density rather than the
fluctuation, and an explicit time integrator adds a small increment to it at every stage. As the time step is refined, these increments
become comparable to the spacing, and each stored value is then rounded by an amount that depends
on the increment rather than on chance. The rounding errors therefore do not average out, and they
accumulate in proportion to the number of steps. In the model problem of Section~\ref{sec:osc},
which has exactly these magnitudes, the stored density does not change at all on $82\%$ of the
steps at $\Delta t=10^{-4}$.

\begin{figure}[!t]\centering
\includegraphics[width=\textwidth]{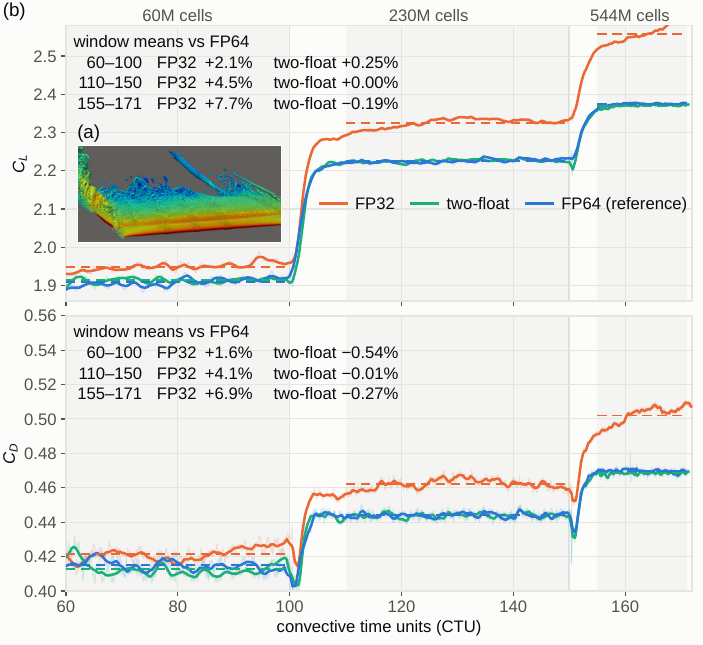}
\caption{Wall-modeled large-eddy simulation of the high-lift Common Research Model (CRM-HLS) of the
Sixth AIAA High-Lift Prediction Workshop (HLPW-6) at $30^\circ$ angle of attack and Mach $0.2$
(Section~\ref{sec:hlpw}). (a) Inset: isosurface of the $Q$-criterion at $10^{5}$\,s$^{-2}$, colored
by velocity magnitude from $0$ (blue) to $200$\,m\,s$^{-1}$ (red), for the two-float build on the
$544$M-cell grid; the FP64 and FP32 fields look the same (Fig.~\ref{fig:edtrio}). (b) Lift (top) and
drag (bottom) histories in convective time units (CTU) for the FP32 (orange), FP64 (blue) and
two-float (green) builds, which differ in how the state is stored between stages. Vertical lines
mark the grid refinements, each halving the finest cell size and the time step; shaded bands are the
averaging windows, dashed lines their means, and the annotations the difference of each mean from
FP64 in percent. Light lines are instantaneous values and heavy lines moving averages}
\label{fig:trio}
\end{figure}

While the opportunity is clear, an important question is to what extent single precision can be
used without sacrificing the accuracy that unsteady simulations require. Several studies have
found the effect to be small. Homann et al.\ \cite{homann2007} examined the floating-point
precision of pseudo-spectral direct numerical simulations of turbulence, Witherden and Jameson
\cite{witherden2020} used FP32 for high-order implicit large-eddy simulations with the flux
reconstruction method, and V\'a\v{n}a et al.\ \cite{vana2017} ran an operational weather
prediction model in single precision; in each case single precision was found adequate for the
quantities of interest. Brogi et al.\ \cite{brogi2024} measured the effect of single and mixed
precision on accuracy and performance in incompressible, compressible and multiphase OpenFOAM
solvers. Most recently, Karp et al.\ \cite{karp2025} compared four scale-resolving
solvers on four flows, with formats emulated down to eight bits, and concluded that IEEE single
precision can be used for the entire simulation, although rounding the stored state proved more
damaging than rounding individual terms of the equations. Sikl\'osi et al.\ \cite{siklosi2025}
observed in an explicit finite-difference solver for compressible flow that, with every array in the
same reduced precision, simulations of the Taylor--Green vortex were slightly less accurate at
smaller time steps, whereas with the state held in higher precision the time step had no effect. The
analysis presented here accounts for these observations, and Section~\ref{sec:earlier} explains why
the earlier findings are compatible with the results reported below.

This work began with a production campaign. In wall-modeled large-eddy simulations (WMLES) of a
high-lift wing, FP32 runs consistently overpredicted the lift relative to FP64, by an amount that
grew with each grid refinement (Fig.~\ref{fig:trio}). Model problems localized the bias to the store
of the Runge--Kutta stages, and the analysis below explains it: the rounding of the stored state, and
only that rounding, grows in effect as the time step is refined. Past a turnaround step, refinement
therefore increases the error, in FP64 as well as in FP32, so that precision and time step cannot be
chosen independently.

The main contributions of this work are as follows.
\begin{itemize}
\item \textbf{Analysis.} We formulate a modified-equation model for Runge--Kutta methods in
Shu--Osher form and for linear multistep methods, building on the classical result that rounding
error can grow as $1/\Delta t$ \cite{henrici1962} and on observations of that growth in geophysical
models \cite{harvey2016,klower2022}. It separates the rounding of the stored state, the only term
with a negative power of $\Delta t$, from the rounding of the right-hand side, which sets a floor.
We bound the state rounding for strong-stability-preserving schemes and identify two mechanisms
that keep its sign from changing, one of which extends the rounded-coefficient error of Hairer et
al.\ \cite{hairer2008} to coefficients that multiply the state. We also give a scale rule that
indicates which variables are affected, and check the analysis exactly on a conserved total.
\item \textbf{Method.} We store each state variable as the unevaluated sum of two FP32 numbers
\cite{dekker1971}, updated with FP32 error-free transformations \cite{muller2018,joldes2017}, while
every right-hand-side evaluation, boundary condition, communication and input/output operation
reads the same FP32 array as before. Unlike the mixed-precision mode of OpenMM \cite{eastman2017},
which updates a similar pair in FP64 arithmetic, it requires no FP64 arithmetic. It needs eight
bytes per stored variable rather than the twelve of an FP64 state kept alongside its FP32 image.
Unlike a compensation term, it is closed under the scalar multiplications of the Shu--Osher form.
We give a storage rule for general Runge--Kutta and multistep schemes and the steps needed to
retrofit it into an existing production solver.
\item \textbf{Generality.} Four model problems from different fields show that the behavior is not
specific to fluid dynamics. Among the remedies that change only the time integrator, only the
two-float state reaches the rounding floor in the usual form of each of SSPRK(3,3), classical RK4,
low-storage SSPRK(10,4) and BDF2; compensation does so only in increment form.
\item \textbf{Production scale.} In WMLES of the Sixth AIAA High-Lift Prediction Workshop wing,
single precision biases the integrated loads by an amount that roughly doubles with each grid
refinement. The window means of the two-float state agree with the FP64 reference to within
$0.6\%$, and the solver runs $1.6$ to $2.7$ times faster than FP64 on a GPU with half-rate FP64.
\end{itemize}
For scale-resolving simulations of aircraft, whose integrated loads are reported to one part in
$10^{4}$, this removes a bias that grows with grid refinement while keeping the speed and memory
savings of single precision. The same holds for other unsteady simulations of states with a large
mean, such as weather, climate and ocean models.

The remainder of this section reviews related work. Section~\ref{sec:modeq} analyzes rounding in
stage-based integrators, Section~\ref{sec:method} describes the two-float state,
Section~\ref{sec:toys} presents the model problems and alternative remedies,
Section~\ref{sec:production} the production results and their cost, and
Section~\ref{sec:discussion} relates the results to earlier studies, gives recommendations and
discusses limitations, before the conclusions in Section~\ref{sec:conclusions}.

\subsection{Related work}\label{sec:prior}

That rounding error in a one-step method can grow as $1/\Delta t$ under refinement is a classical
result \cite{henrici1962}; it occurs in FP64 as well, where it is rarely noticeable
\cite{croci2026}. Henrici also modeled the individual rounding errors statistically: a
systematic part accumulates in proportion to the number of steps, whereas a random part with
independent, zero-mean errors accumulates as a random walk, in proportion to $\Delta t^{-1/2}$, and
Hull and Swenson \cite{hull1966} tested such models experimentally. The classical bound is a worst
case and does not by itself establish that the growth occurs. Harvey and Verseghy
\cite{harvey2016} found in a land-surface model that, at a given precision, smaller time steps
increase the rounding error of the stored temperature even as they reduce the truncation error, and
Kl\"ower et al.\ \cite{klower2022} observed that the error of a 16-bit shallow-water model grows as
the time step is refined. In long orbital integrations in FP64, the error likewise has an optimum
step below which it rises \cite{rein2015whfast}. Croci and Giles
\cite{croci2023} showed for the heat equation that round-to-nearest can make the solution stagnate
and the rounding error grow as $\varepsilon/\Delta t$ at a fixed horizon, whereas stochastic rounding
avoids both.

Most work on mixed-precision time integration lowers the precision of the evaluations while the
state is held in FP64. Grant \cite{grant2022} introduced perturbed Runge--Kutta methods, in which
expensive stages are computed in low precision and corrected with high-precision ones. Burnett et al.\ \cite{burnett2024} analyzed the stability of such additive mixed-precision
methods, and Field et al.\ \cite{field2023} used mixed precision in a GPU-accelerated WENO method.
Croci and
Rosilho de Souza \cite{croci2022} computed the stability-only stages of explicit stabilized
Runge--Kutta methods in low precision, and Balos et al.\ \cite{balos2023} did the same for the
Krylov approximations used in exponential integrators. Kelley \cite{kelley2022} analyzed Newton's
method in mixed precision, in which inner iterations are carried out in single precision under a
double-precision correction. Andena and Croci \cite{croci2026} recently surveyed reduced- and
mixed-precision Runge--Kutta methods, including compensated summation and stochastic rounding, and
noted that holding the state and adding the update in high precision, with the update computed in
low precision, bounds the rounding error by $O(u+u_hT/\Delta t)$, where $u$ and $u_h$ are the unit
roundoffs of the two formats. In
each of these designs the state, and therefore the error term isolated here, is held in the higher
of two precisions, which in practice is FP64. The designs closest to the present one also hold the
state more precisely than the evaluations. Kl\"ower et al.\ \cite{klower2020} kept the prognostic
variables of a shallow-water model in Float32 while computing the tendencies in Float16 or
BFloat16. The mixed-precision mode of the molecular-dynamics code OpenMM \cite{eastman2017} computes
forces in FP32 and stores each position as an FP32 value with an FP32 correction, updated in FP64
arithmetic. GPU N-body codes used two-float arithmetic for predicted positions and position
differences, while the corrector remained in FP64 \cite{gaburov2009,nitadori2012}. The present
approach updates its pair with FP32 error-free transformations only, at every stage of Runge--Kutta
and multistep schemes; in the terms of the bound above, it gives $u_h$ of order $u^2$ without a
second hardware format.

Protecting the update of the state is itself an established idea. Gill \cite{gill1951} proposed a
Runge--Kutta process for fixed-point machines with this aim, and Kahan \cite{kahan1965} and
M{\o}ller \cite{moller1965} introduced compensated summation \cite[Ch.~4]{higham2002}. Compensated
updates are standard practice in long orbital integrations \cite[Sec.~VIII.5]{hairer2006},
\cite{quinn1990,hairer2008,rein2015}. More
recently, Kl\"ower et al.\ \cite{klower2022} ran a 16-bit shallow-water model with a compensated
prognostic update, Lai et al.\ \cite{lai2025} applied the same compensation to the FP32 state
update of the MPAS-A atmospheric dynamical core, and Dawson et al.\ \cite{dawson2018} retained
higher precision in the part of a land-surface model that accumulates small increments. In
molecular dynamics, Colberg and H\"ofling \cite{colberg2011} traced spurious glassy dynamics to
single-precision rounding of the particle positions and corrected it with a compensating
low-order part in the position update, and Le Grand et al.\ \cite{legrand2013} used fixed-point
accumulation for similar reasons. All of these works protect updates in increment form, the form in which compensated summation is
also formulated for Runge--Kutta methods in \cite{croci2026}, in which
no scalar multiplies the state, and most do so with a single compensation term. As Section~\ref{sec:progression}
explains, this is not sufficient for the Shu--Osher form used by
strong-stability-preserving solvers, in which the significant rounding errors occur in scalar
multiplications of the state before any addition that a compensation term could protect. The same difficulty arises in multistep
methods for orbital mechanics, in which a coefficient multiplying the state amplifies the rounding
and compensating only the final sum has little effect; Grazier et al.\ \cite{grazier2005} and
Console and Hairer \cite{console2014} removed it by rewriting the method in a summed or difference
form. The
two-float representation used here is closed under those multiplications because it represents a
value rather than a correction.

Two-float arithmetic dates to Dekker \cite{dekker1971}. Rigorous error bounds for its building
blocks were given by Jolde\c{s} et al.\ \cite{joldes2017}, Muller and Rideau
\cite{mullerrideau2022} and Boldo et al.\ \cite{boldo2017}, and it was implemented on GPUs as soon
as they became programmable \cite{thall2006,dagraca2006}. For PDE solvers on GPUs, however, the
earlier assessment was unfavorable: G\"oddeke et al.\ \cite{goddeke2007} used emulated double
precision throughout finite-element solvers and found it uncompetitive with mixed-precision
iterative refinement. The present work does not revisit that comparison. It instead applies
two-float arithmetic only where the state is stored, leaving every flux, closure and boundary
evaluation in ordinary FP32, and the resulting cost is a small fraction of that of emulation
throughout (Section~\ref{sec:cost}). FP64 emulation on low-precision matrix units
\cite{ozaki2012,ootomo2022} is complementary to this approach, since it accelerates dense matrix
products rather than the matrix-free flux evaluations of an explicit solver.

Storing the state more precisely than the derivative evaluations is therefore an old idea
\cite{gill1951}. What is new here is narrower: the identification of the store rounding as the only
term that grows under refinement, and of the coefficient rounding, known in Butcher form
\cite{hairer2008} and for multistep methods \cite{console2014}, as a mechanism that keeps its sign
when the coefficients multiply the state; a state held as two FP32 numbers without FP64 arithmetic;
and the demonstration of both the problem and the remedy at production scale.

\section{Analysis of rounding in stage-based time integrators}\label{sec:modeq}

\subsection{Preliminaries}\label{sec:acc}

Let $\fl(\cdot)$ denote rounding to the nearest FP32 number, and let $\ulp(x)$, the unit in the
last place, denote the spacing of FP32 numbers at $x$,
$\ulp(x)=2^{-23}\,2^{\lfloor\log_2|x|\rfloor}$, which lies between $0.6\times10^{-7}|x|$ and
$1.2\times10^{-7}|x|$ depending on where $x$ falls between consecutive powers of two. The FP32 unit roundoff is
$\varepsilon=2^{-24}$, half an ulp at magnitude one. A correctly rounded store introduces an error
of at most $\tfrac12\ulp$ of the stored value, and an increment smaller than that cannot change
the stored value.

Three simple experiments with a running total illustrate the mechanism studied in this paper.
Near $10^{6}$ the spacing of FP32 numbers is $0.0625$. Adding $0.03$ to $10^{6}$ a thousand times
leaves the total at $10^{6}$, because each increment is smaller than half the spacing and is lost;
the exact result is $1{,}000{,}030$. Adding $0.05$ instead, which is slightly more than half the
spacing, rounds every increment up to $0.0625$, an error of $+25\%$ with the same sign each time,
and after $10^{4}$ additions the total is $1{,}000{,}625$ rather than $1{,}000{,}500$. Halving the
increment and doubling the number of additions, which is analogous to halving $\Delta t$,
increases the accumulated error. Finally, if the total is kept as an FP32 number plus a small FP32
correction, with each addition carried out in extended precision and the result split on store,
the exact result is recovered from the same FP32 increments. These experiments show that the
rounding error is set by the magnitude of the stored value rather than by the increment, and that
only the rounding of the running total accumulates: an error in an increment is used once and
discarded, whereas an error in the total persists. A time integrator performs the same operation
when it adds a small increment to a large stored state at every stage, for every unknown. The remainder
of this section makes this observation precise.

\subsection{Forward Euler}\label{sec:euler}

After spatial discretization, an evolution problem takes the form
\begin{equation}
\frac{du}{dt}=L(u),
\label{eq:mol}
\end{equation}
which we consider componentwise in the state $u$. The mechanism is simplest for the forward
Euler method. A computed step can be written as
\begin{equation}
\hat u^{n+1}=\hat u^{n}+\Delta t\,\big[L(\hat u^{n})+\lambda_n\big]+\delta_n ,
\label{eq:euler}
\end{equation}
where $\delta_n$, at most half an ulp of the stored value, is the rounding error of the stored
state, and $\lambda_n$ is the rounding error of the right-hand side (RHS). Over a horizon $T$, after
$N=T/\Delta t$ steps, the state rounding errors sum to $T\bar\delta/\Delta t$, where $\bar\delta$
is their mean, the RHS rounding errors sum to $T\bar\lambda$, and the truncation error is
$C\Delta t$. Refining the step reduces the truncation error and leaves $T\bar\lambda$ unchanged,
but increases the accumulated state rounding error. The total error therefore reaches a minimum
near $\Delta t=(T\bar\delta/C)^{1/2}$ and then increases as $1/\Delta t$. This is Henrici's
rounding bound \cite{henrici1962}. Whether the bound is attained depends on the mean
$\bar\delta$, which is discussed in Section~\ref{sec:sign}. Runge--Kutta methods store several
stage states per step, and the analysis below accounts for them.

\subsection{Shu--Osher form and path weights}\label{sec:so}

An explicit Runge--Kutta method in Shu--Osher form \cite{shu1988} advances Eq.~\eqref{eq:mol} by
\begin{equation}
u^{(i)}=\sum_{k<i}\Big[\alpha_{ik}\,u^{(k)}+\Delta t\,\beta_{ik}\,L\big(u^{(k)}\big)\Big],
\qquad \textstyle\sum_k \alpha_{ik}=1,
\label{eq:so}
\end{equation}
with $u^{(0)}=u^{n}$ and $u^{n+1}=u^{(s)}$. For strong-stability-preserving (SSP) schemes every
row is convex, $\alpha_{ik}\ge0$ \cite{gottlieb2001}. The scheme used in the production solver,
SSPRK(3,3), has rows $(1)$, $(\tfrac34,\tfrac14)$ and $(\tfrac13,\tfrac23)$, with each $\beta$
equal to the last $\alpha$ of its row. In contrast to the Butcher-tableau form, in which
increments are accumulated, this form multiplies the state itself by scalars and stores more than
one stage state at the magnitude of the state. Both properties are important in
Section~\ref{sec:method}.

Two kinds of error separate the computed stages $\hat u^{(i)}$ from their exact counterparts. The
state rounding error $\delta_i$ includes everything rounded away in forming and storing stage
$i$. It is of order $\ulp(\bar u)$, where $\bar u$ is the stored magnitude: half an ulp for the
store itself, and up to half an ulp more for each product and sum formed at the magnitude of the
state. The RHS rounding error $\lambda_k$ is the error of the FP32 evaluation of $L$ at stage $k$,
which is dominated by the quantization of the state that $L$ reads. The computed stages satisfy,
exactly,
\begin{equation}
\hat u^{(i)}=\sum_{k<i}\Big[\alpha_{ik}\,\hat u^{(k)}
 +\Delta t\,\beta_{ik}\,\big(L(\hat u^{(k)})+\lambda_k\big)\Big]+\delta_i .
\label{eq:pert}
\end{equation}
Subtracting Eq.~\eqref{eq:so} from Eq.~\eqref{eq:pert} gives a recursion for the stage errors
$e^{(i)}=\hat u^{(i)}-u^{(i)}$ within one step, starting from $e^{(0)}=0$,
\begin{equation}
e^{(i)}=\sum_{k<i}\alpha_{ik}\,e^{(k)}+\delta_i+\Delta t\sum_{k<i}\beta_{ik}\Big[\lambda_k+L(\hat
u^{(k)})-L(u^{(k)})\Big],
\label{eq:errrec}
\end{equation}
in which the bracket contains the RHS rounding error and the change in $L$ caused by the stage
error, which is proportional to $e^{(k)}$ for smooth $L$. Every term involving $L$ carries a factor
of $\Delta t$, so to leading order the error of each stage is a convex combination of the errors of
its predecessors plus its own state rounding error. Unrolling Eq.~\eqref{eq:errrec} backward from
the last stage, $\delta_i$ reaches $u^{n+1}$ along every chain of stages $i\to j\to\dots\to s$,
multiplied by the product of the coefficients $\alpha$ along the chain. The sum over all chains is
the path weight $w_i$, given by the backward recursion
\begin{equation}
w_s=1,\qquad w_i=\sum_{j>i}\alpha_{ji}\,w_j .
\label{eq:weights}
\end{equation}
For convex rows, the coefficients of each row form a probability distribution over the preceding
stages, and Eq.~\eqref{eq:weights} gives the probability $w_i$ that a backward chain starting at
the final stage passes through stage $i$. Hence $0\le w_i\le1$, so that no SSP scheme amplifies a
state rounding error, although the error of a given stage may be attenuated. For SSPRK(3,3),
$w_3=1$, $w_2=\tfrac23$ and $w_1=\tfrac14\cdot\tfrac23=\tfrac16$, and the weight of $u^{n}$ is
$\tfrac16+\tfrac12+\tfrac13=1$, as consistency requires. At the next order, $\lambda_k$ enters stage
$j$ with weight $\Delta t\,\beta_{jk}$ and is carried from there with weight $w_j$, so that its
total weight is $b_k=\sum_{j>k}\beta_{jk}w_j$. For SSPRK(3,3) these weights are
$(\tfrac16,\tfrac16,\tfrac23)$, and for any consistent scheme they sum to one. The change in $L$
caused by a stage error is smaller than the state rounding error that caused it by a factor of
order $\Delta t\,|\partial L/\partial u|$. Near the stability limit this factor is of order one,
and the path weights become polynomials in $\Delta t\,\partial L/\partial u$; however, the factor
vanishes under refinement, and the regime in which the state rounding error dominates is precisely
the one in which it is small. To leading order, the error added in one step is therefore
\begin{equation}
e^{(s)}=\sum_i w_i\,\delta_i\;+\;\Delta t\sum_k b_k\,\lambda_k\;+\;\text{higher-order terms}.
\label{eq:inject}
\end{equation}

\subsection{Modified equation}\label{sec:modeqn}

Each step adds $e^{(s)}$ to the result that exact arithmetic would give. Since a forcing $f$ in
$du/dt=L(u)+f$ changes the state by $\Delta t\,f$ per step, adding $e^{(s)}$ at every step is
equivalent to forcing the equation by $e^{(s)}/\Delta t$. By backward error analysis
\cite{hairer2006}, the truncation error of the scheme is a forcing of the same kind. Dividing
Eq.~\eqref{eq:inject} by $\Delta t$, which is the classical local-error-per-unit-step argument, the
computed trajectory follows, to leading order, the modified equation
\begin{equation}
\frac{du}{dt}=L(u)
\;+\;\underbrace{C(u)\,\Delta t^{\,p}}_{\tau\text{, truncation}}
\;+\;\underbrace{\frac{1}{\Delta t}\sum_i w_i\,\delta_i}_{R_u\text{, state rounding}}
\;+\;\underbrace{\sum_k b_k\,\lambda_k}_{R_L\text{, RHS rounding}},
\label{eq:modeq}
\end{equation}
where $p$ is the design order and $C(u)$ the truncation coefficient. The rounding terms are not
smooth functions of $u$, so Eq.~\eqref{eq:modeq} should be regarded as a model rather than as a
backward-error result in the strict sense; a general stability theory for reduced- and
mixed-precision methods applied to nonlinear systems remains open for the same reason, since
rounding errors act as non-smooth perturbations \cite{croci2026}. Section~\ref{sec:theorems} checks one of its predictions
exactly. It describes the mean part of the rounding forcing,
whose accumulation is examined in Section~\ref{sec:sign}, while the zero-mean remainder
accumulates as a random walk. Equation~\eqref{eq:modeq} contains nothing specific to fluid
dynamics: $L$ is an arbitrary right-hand side, and each forcing term is present for any
integrator that stores its state once per stage.

Each branch of a step-size study corresponds to one term of this equation. The truncation term
$\tau$ dominates the branch on which the error decreases with $\Delta t$. The state rounding term
$R_u$ is the only term with a negative power of $\Delta t$, a divergent forcing that grows as
$\Delta t\to0$ until the increments fall below half an ulp and the state stagnates, and it
dominates the branch on which the error increases. The RHS rounding term $R_L$, whose factor of
$\Delta t$ has canceled, is a bounded forcing whose effect does not depend on $\Delta t$; it does
depend on the problem and the horizon, and for the oscillator of Section~\ref{sec:osc}, whose
rounding perturbation is conservative and averages out over a cycle, it is flat. It sets the floor
reached by any method that does not have the $R_u$ term. The magnitude of $R_u$ is set
by the stored magnitude rather than by the signal. For SSPRK(3,3), the store rounding alone
contributes at most
$(\tfrac16+\tfrac23+1)\tfrac12\ulp(\bar u)/\Delta t=\tfrac{11}{12}\ulp(\bar u)/\Delta t$. No finite
store eliminates the $\Delta t^{-1}$ term. The aim of the method in Section~\ref{sec:method} is
therefore to make it smaller by a factor of about $1/\varepsilon$, which is $2^{24}$ for FP32, so that it lies below the
floor set by $R_L$ over the steps of interest, while leaving the other terms unchanged.

Balancing the truncation and state rounding terms gives the location of the turnaround. We write
the error at a fixed horizon as $E(\Delta t)=k_\tau\Delta t^{p}+k_R/\Delta t$, where $k_\tau$ is
the global truncation constant and $k_R$, which is proportional to the unit roundoff of the format
for a given rounding constant, represents the integrated state rounding. The error is smallest at
$\Delta t^{*}=(k_R/p\,k_\tau)^{1/(p+1)}$, where $E_{\min}\propto k_R^{\,p/(p+1)}$. The turnaround
step therefore scales with the $(p+1)$-th root of the unit roundoff and the minimum error with its
$p/(p+1)$ power, so that for two formats the ratio of the minimum errors equals the $p$-th power of
the ratio of the turnaround steps. A more precise number format thus moves the turnaround to
smaller steps but does not remove it, and the same holds for a two-float state, whose own store
rounding is of order $\varepsilon^2\bar u$. The difference is that a wider format lowers the floor
$R_L$ together with the $\Delta t^{-1}$ term, whereas a two-float state with an FP32 right-hand
side lowers only the latter, which then lies beneath the floor until the step is many orders of
magnitude smaller. This prediction is tested in Section~\ref{sec:osc}.

\paragraph{Butcher form and linear multistep methods} An integrator written in Butcher-tableau
form, with the state updated by increments, is also affected. It stores the same state at the
same magnitude and rounds it at the same $\tfrac12\ulp(\bar u)$; only the number of roundings that
reach the end of a step with weight one changes, and the bound for the store rounding decreases
from $\tfrac{11}{12}$ to $\tfrac12\ulp(\bar u)/\Delta t$ without changing the dependence on
$\Delta t$. A Butcher-form version of the five-station oscillator of Section~\ref{sec:refstate}
still changes the oscillation amplitude by $+4.7\%$ at $\Delta t$ and by $+51\%$ at $\Delta t/4$.
These errors are smaller than those of the Shu--Osher form by more than the ratio of the two
bounds, because the rounding constant differs between the forms, but they grow under refinement
in the same way. Linear multistep methods have the form of Eq.~\eqref{eq:so} with retained time
levels in place of stages. For example, BDF2 is
$u^{n+1}=\tfrac43u^{n}-\tfrac13u^{n-1}+\tfrac23\Delta t\,L(u^{n+1})$, whose coefficients sum to one
as those of a Shu--Osher row do, so that the derivation applies with one qualification: the rows
are affine but not convex, and the path weights can exceed one.

\subsection{Exposure of individual variables}\label{sec:scales}

For a state with a large mean $\bar u$ and a small fluctuation of amplitude $A$, the forcing $R_u$
is of order $\ulp(\bar u)/\Delta t$. After $N=T/\Delta t$ steps the accumulated error is at most of
order $N\ulp(\bar u)$, or $N\ulp(\bar u)/A$ relative to the fluctuation. The ratio
\begin{equation}
\frac{\ulp(\text{stored mean})}{\text{fluctuation}}
\label{eq:scale}
\end{equation}
therefore indicates in advance which variables are affected first, and multiplied by the number of
stage updates it gives an upper bound on the relative error rather than a prediction of it. We
refer to Eq.~\eqref{eq:scale} as the scale rule. It acts as a prefactor rather than a threshold,
and it identifies the exposed variables without ranking the others: a variable without a large mean is still affected, but at
smaller time steps, as the shallow-water results of Section~\ref{sec:sw} show.

The ratio is dimensionless: rescaling a variable rescales its ulp by the same factor, up to the
factor of two set by where the stored mean falls between consecutive powers of two, so
nondimensionalization cannot remove it. It does depend on the physics of the state. For a
compressible flow at sea level and Mach $0.2$, the density is stored near $1.225$ with fluctuations
of about $10^{-3}$, and the total energy near $2.56\times10^{5}$ with fluctuations of about $600$.
Pressure and temperature, which these variables encode, therefore have ratios near $10^{-4}$,
whereas the ratios for the momenta are about a hundred times smaller (Table~\ref{tab:fivevar}).
The integrated loads on an aerodynamic body are surface integrals of the pressure fluctuation,
which explains why they show the bias reported in Section~\ref{sec:production}.

\subsection{Why the rounding errors do not cancel}\label{sec:sign}

The bound on $R_u$ becomes a growth law only if the state rounding errors keep the same sign, so
that they add linearly with the number of steps rather than cancelling as a random walk. Two
mechanisms keep the sign from changing.

The first is the size of the increment. When a stage changes the stored value by only a few ulp,
the rounding error is determined by the increment itself rather than by chance, and the same sign
recurs. A resolved flow changes slowly from step to step, so its increments, and therefore the
rounding errors, repeat. The mean then drifts by $\gamma\,\ulp(\bar u)$ per stage, with $\gamma$
between $-\tfrac12$ and $\tfrac12$ determined by the distribution of the increments, and by
$\gamma N\ulp(\bar u)$ after $N$ stages. Since $\gamma$ depends on the shape of that distribution
rather than on its size, the sign of the bias depends on the flow and cannot be calibrated out.
Croci and Giles \cite{croci2023} analyzed related behavior, stagnation under round-to-nearest, for
the heat equation.

The second mechanism is the rounding of the coefficients of the scheme. The FP32 representations of
$\tfrac13$ and $\tfrac23$ in the final row of SSPRK(3,3) sum to $1+2^{-25}$. This coefficient defect
amounts to a gain applied to the whole state once per step; for a state of magnitude $\bar u$ this is a systematic
change of between a quarter and half an ulp per step, depending on where $\bar u$ falls between
consecutive powers of two. Hairer et al.\ \cite{hairer2008} identified the corresponding systematic error of rounded
Runge--Kutta coefficients in long-time integrations, where the coefficients multiply $\Delta t\,f$,
and avoided it by representing each coefficient as the sum of two floating-point numbers;
Console and Hairer \cite{console2014} gave the corresponding advice for multistep methods,
scaling rational coefficients to integers. In the
Shu--Osher form the defect multiplies the state itself, so that its contribution at a fixed final
time grows as $\Delta t^{-1}$. Coefficients
that multiply increments, namely the $\beta_{ik}$ and $\Delta t$, also produce a systematic error,
the one removed by Hairer et al., but it does not grow under refinement: its contribution per step
is proportional to $\Delta t$, so that it accumulates to an amount independent of $\Delta t$ at a
fixed final time, and it is smaller still because it multiplies a quantity of the size of the
fluctuation rather than the mean.

Which mechanism dominates depends on the size of the increment. When the increments are many ulp,
the coefficient defect is the only systematic error; the model problem of
Section~\ref{sec:pendulum} is an example. When they are close to one ulp, as for the compressible state, the state rounding error
dominates, and renormalizing the row so that it sums to one exactly in FP32 does not help
(Section~\ref{sec:coeffs}). An intermediate case occurs when the sign is maintained over a run of
steps and then reverses. If the fractional position of the increment on the rounding grid drifts
by $a\omega^{2}\Delta t^{2}$ per step for a wave of amplitude $a$ and frequency $\omega$, a run
lasts about $\ulp/(a\omega^{2}\Delta t^{2})$ steps, and runs of that length accumulating as a
random walk give an error at a fixed horizon that grows as $\Delta t^{-3/2}$. The shallow-water
results of Section~\ref{sec:sw} are consistent with this exponent. It is transitional: when a run
lasts longer than the whole integration the sign is fixed and the growth reverts to $\Delta t^{-1}$,
and when runs are short it relaxes toward the random-walk value $\Delta t^{-1/2}$; the fitted
shallow-water exponents of $-0.8$ and $-1.5$ lie between these limits. In each case the error grows as
the step is refined, and the two-float state of Section~\ref{sec:method} reduces the growing term
by a factor of about $2^{24}$.

\subsection{A check on a conserved total}\label{sec:theorems}

The modified equation \eqref{eq:modeq} is a model, but one of its predictions can be checked
exactly. Many discretizations conserve a total, such as the mass summed over the cells of a
finite-volume scheme on a periodic or closed domain. In exact arithmetic that total does not change
at all, so any change in the computed total is caused by rounding alone, with no truncation error to
separate from it. Summing Eq.~\eqref{eq:pert} over the cells removes the right-hand side, because
the fluxes cancel in pairs, and leaves only the rounding errors of each stage, carried to the end of
the step with the path weights of Eq.~\eqref{eq:weights}. The change of the total $Q$ over a run of
$N=T/\Delta t$ steps is therefore exactly
\begin{equation}
Q_N-Q_0=\frac{T}{\Delta t}\,\bar\sigma+T\bar\mu ,
\label{eq:driftsum}
\end{equation}
where $\bar\sigma$ is the average per step of the weighted state rounding errors, summed over the
cells, and $\bar\mu$ the corresponding average of the right-hand-side rounding errors. Unlike
Eq.~\eqref{eq:modeq}, this is an identity rather than a model. It shows that the state rounding is
the only part multiplied by $1/\Delta t$: halving the step doubles its contribution unless its
average is zero. Part of that average is never zero. The FP32 coefficients of the last row of
SSPRK(3,3) sum to $1+2^{-25}$, which changes the total by $2^{-25}Q$ at every step, always in the
same direction; the two-float coefficients of Section~\ref{sec:coeffcarrier} reduce this to
$2^{-49}$.

As a check, consider periodic linear advection with first-order upwind fluxes on $64$ cells, with a
mean of $1.225$ and a sine wave of amplitude $10^{-3}$, integrated with SSPRK(3,3) to $T=1$.
Equation~\eqref{eq:driftsum} reproduces the computed FP32 change of the total exactly. The change
grows from $5.4\times10^{-4}$ to $1.4\times10^{-2}$ as $\Delta t$ is reduced from $4\times10^{-3}$
to $2.5\times10^{-4}$, and the coefficient defect alone, $2^{-25}Q\approx2.3\times10^{-6}$ per
step, is close to the measured $2.1\times10^{-6}$ per step at the largest step. With the two-float
state, modelled by evaluating each row in FP64 and splitting the result, the change stays below
$2\times10^{-11}$.

The same argument, combined with the standard stability estimate for one-step methods
\cite{henrici1962}, bounds the error of the solution itself by the three terms of
Eq.~\eqref{eq:modeq}: a truncation term, a right-hand-side term that does not depend on $\Delta t$,
and a state rounding term proportional to the rounding per stage divided by $\Delta t$. Only the last
depends on how the state is stored, and the two-float state of Section~\ref{sec:method} makes it
about five orders of magnitude smaller (Section~\ref{sec:bound}). Neither result shows that the average
rounding error is nonzero for every problem: Section~\ref{sec:sign} gives the mechanisms that keep
it away from zero, Sections~\ref{sec:toys} and~\ref{sec:production} measure their effect, and a
general proof remains open \cite{croci2026}.

\section{Two-float state representation}\label{sec:method}

Equation~\eqref{eq:modeq} indicates what should be changed and what can be left as it is. The
state rounding error is the only term that accumulates, so only the stored state requires higher
precision. The analysis does not specify how that precision is provided. An FP64 state provides it
where FP64 arithmetic is fast and memory is not a constraint. The approach used here instead
builds the state from two numbers of the working format, so that it remains effective where FP64
is slow or unavailable, requires the same eight bytes per stored variable as an FP64 state, and
keeps the existing FP32 array as its leading part, so that every part of the code that reads the
state is unchanged.

\subsection{Representation and update}\label{sec:rep}

Each state variable $u$ that persists across a stage boundary is stored as an unevaluated
pair of FP32 numbers \cite{dekker1971},
\begin{equation}
u=\hi+\lo,\qquad \hi=\fl(\hi+\lo),\qquad |\lo|\le\tfrac12\ulp(\hi),
\label{eq:pair}
\end{equation}
so that $\hi$ carries the leading seven significant digits and $\lo$ the next seven, about $48$
significant bits in total. This is less than the $53$ bits of FP64, but much more than the $24$
bits of FP32. We therefore do not claim double precision, but rather agreement with the FP64
reference in each quantity, as measured below. The $\hi$ float is the existing state array: the
flux kernels, boundary conditions, wall model, MPI exchanges and file input and output read the
same FP32 state as before, and $\lo$ holds what a single FP32 store would have discarded.

\begin{figure}[!tbp]\centering
\includegraphics[width=\textwidth]{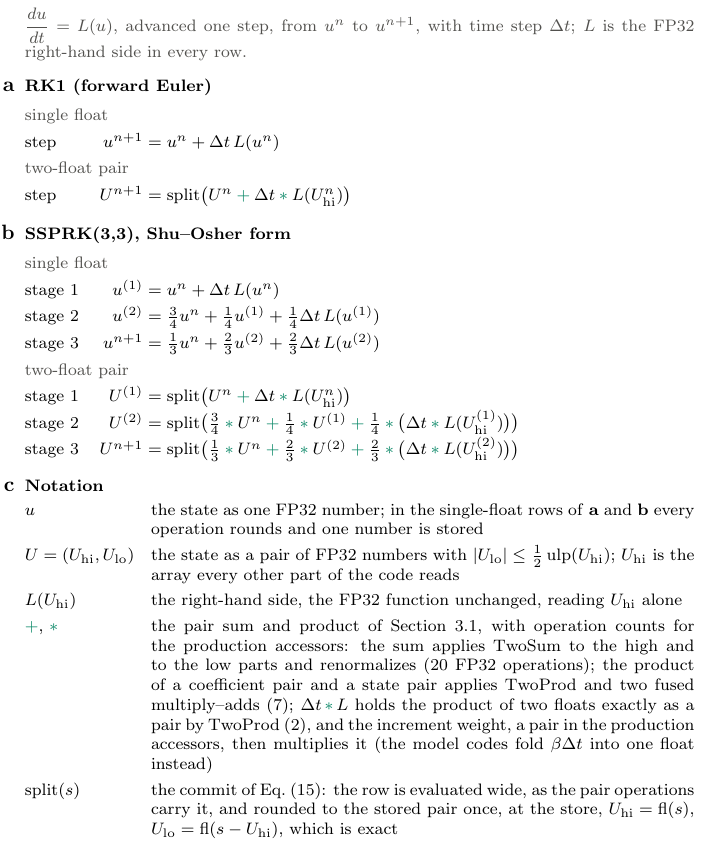}
\caption{FP32 and two-float state updates for (a) one step of forward Euler and (b) the three stages
of SSPRK(3,3) in Shu--Osher form, one line per stage. In both, the right-hand side $L$ is the same
FP32 function and reads only the leading float $U_{\hi}$. (c) Notation: the pair
$U=(U_{\hi},U_{\lo})$; the two-float sum and product of Eqs.~\eqref{eq:twosum}--\eqref{eq:row},
drawn in green, with their FP32 operation counts; and the split of Eq.~\eqref{eq:split}}
\label{fig:carrier}
\end{figure}

All updates of the state are made through two accessor functions built from the classical
error-free transformations \cite{muller2018,joldes2017}, which compute a rounded FP32 result and
then recover, with a few more FP32 operations, the exact error of the rounding. With $\oplus$,
$\ominus$ and $\otimes$ denoting correctly rounded FP32 operations and $\mathrm{fma}$ the fused
multiply--add,
\begin{align}
(s,e)=\mathrm{TwoSum}(a,b):&\quad s=a\oplus b,\;\; v=s\ominus a,\;\;
  e=\big(a\ominus(s\ominus v)\big)\oplus(b\ominus v);\label{eq:twosum}\\
(p,e)=\mathrm{TwoProd}(a,b):&\quad p=a\otimes b,\qquad e=\mathrm{fma}(a,b,-p),\label{eq:tpf}
\end{align}
and each satisfies an exact identity, $s+e=a+b$ and $p+e=ab$. TwoSum does not require $|a|\ge|b|$,
which is important because the momenta pass through zero; the less expensive Fast2Sum, which does
require it, is used only within the published two-float algorithms whose proofs allow it
\cite{joldes2017,mullerrideau2022,boldo2017}. Combining these into the two-float sum and the
products of a pair with a pair and with an FP32 number \cite{joldes2017}, the accessors evaluate
one Shu--Osher row with the state held in two-float form,
\begin{equation}
u^{\mathrm{new}}=\sum_k\alpha_k\,(\hi_k+\lo_k)+\Delta t\sum_k\beta_k\,L_k,
\qquad L_k=L(\hi_k)\ \text{in FP32},
\label{eq:row}
\end{equation}
and store the result with the normalizing split
\begin{equation}
\hi^{\mathrm{new}}=\fl(u^{\mathrm{new}}),\qquad
\lo^{\mathrm{new}}=\fl\big(u^{\mathrm{new}}-\hi^{\mathrm{new}}\big),
\label{eq:split}
\end{equation}
which is again exact. The rounding error that an FP32 store would discard is retained in $\lo$,
while the RHS evaluations $L_k$ remain ordinary FP32 computations on the $\hi$ floats.
Figure~\ref{fig:carrier} compares the FP32 and two-float updates for forward Euler and for
SSPRK(3,3). In the production solver the last two-float addition of a row already returns a
normalized pair, so that the split requires no additional operations. The two-float addition costs
$20$ FP32 operations, and the product of a coefficient pair with a state pair costs $7$: a TwoProd
of the leading parts, two fused multiply--adds for the cross terms, and a Fast2Sum to renormalize.

\paragraph{Precision of the coefficients}\label{sec:coeffcarrier} The Shu--Osher coefficients $\alpha_k$ are stored as two-float constants at compile time, so that
the row-sum defect of their FP32 representations, $1+2^{-25}$ for the final row of SSPRK(3,3), is
reduced to $2^{-49}$; general coefficient tables are renormalized with a two-float reciprocal. The
increment weights $\beta_k$ may be stored as single floats, as in the model codes, or as pairs, as
in the production implementation. The time step is kept as a single float, because its product
with $L_k$ is computed exactly by the fused multiply--add of Eq.~\eqref{eq:tpf} and its error
multiplies a fluctuation rather than the mean. Only the accumulated time $t=\sum\Delta t$ requires
higher precision, for bookkeeping.

\subsection{Relation to the FP32 and compensated updates}\label{sec:progression}

The two-float state can be viewed as the last of three successive ways of storing the result of a
stage, and it is instructive to compare them directly. In the plain FP32 update, every operation of
Eq.~\eqref{eq:so} is rounded,
\begin{equation}
\hat u^{(i)}=\fl\Big(\textstyle\sum_k \fl\big(\alpha_{ik}\hat u^{(k)}\big)
  +\fl\big(\Delta t\,\beta_{ik}L_k\big)\Big),
\label{eq:fp32row}
\end{equation}
and the rounding errors of the products $\alpha_{ik}\hat u^{(k)}$ and of the final store, all at
the magnitude of the state, form the $\delta_i$ of Section~\ref{sec:so}. For the oscillator of
Section~\ref{sec:osc} this update loses $7.2\%$ of the amplitude in twenty periods, and its error
grows as $\Delta t^{-1}$ under refinement.

The established improvement is compensated summation \cite{kahan1965,moller1965}, as used by
Kl\"ower et al.\ \cite{klower2022} and Lai et al.\ \cite{lai2025}. For an update in increment
form, $u^{n+1}=u^{n}+\Delta u$, a second float $c$ carries the rounding error of each addition into
the next one,
\begin{equation}
y=\Delta u\ominus c,\qquad s=u^{n}\oplus y,\qquad c\leftarrow(s\ominus u^{n})\ominus y,
\qquad u^{n+1}=s .
\label{eq:kahan}
\end{equation}
With classical RK4, whose combination rows are all ones, or with SSPRK(3,3) rewritten in increment
form, it gives the same result as the two-float state when the increments span many ulp
(Section~\ref{sec:pendulum}). It is insufficient in two situations. First, in Shu--Osher form the
significant rounding errors occur where scalars multiply the state, in the rounded coefficients
$\alpha_{ik}$ and the products $\alpha_{ik}\hat u^{(k)}$ of Eq.~\eqref{eq:fp32row}, before any
addition that Eq.~\eqref{eq:kahan} could protect: for the oscillator at a mean of $10^{6}$ the
compensated amplitude error is $-1.2\%$, against $-7.2\%$ for FP32; on every step-size sweep of
Fig.~\ref{fig:domains} the compensated error grows with the FP32 error under refinement
(Sections~\ref{sec:sw}--\ref{sec:nbody}). Second, the compensation
term is not part of the stored value, so the stage arguments are formed from the FP32 value alone
and stagnate when the increments are close to one ulp: in increment form the oscillator at a mean
of $10^{6}$ then has an amplitude error of $+13\%$, and in the sub-ulp regime of Fig.~\ref{fig:mech}
the compensated error is erratic and $10$ to $40$ times the two-float floor (Fig.~\ref{fig:mech}b,
purple diamonds). Folding the compensation into every stage value as well, which we call
compensated stages, removes the second shortcoming but not the first: for a scheme written in
increment form it comes within a factor of two of the two-float floor in every model problem (open
triangles in Figs.~\ref{fig:mech}b and~\ref{fig:domains}). At about four FP32 operations per
addition, against twenty for the two-float sum, it is the cheaper choice for an integrator already
in increment form whose state is modified only by the update, but it requires that form, whereas the
two-float state applies unchanged to the Shu--Osher, low-storage and multistep forms in which scalars
multiply the state.

Both shortcomings are removed by treating $u$ and its correction as a single value. The second
float then enters every stage argument through the sum $\hi+\lo$ of Eq.~\eqref{eq:row}, and the
products with the coefficients are computed with error-free transformations, so that the
representation is closed under the multiplications of the Shu--Osher form. This is the two-float
state of Section~\ref{sec:rep}. For the oscillator it reduces the amplitude error to
$1.5\times10^{-4}$, and in every model problem of Section~\ref{sec:toys} it follows FP64 down to
the right-hand-side floor. A further
practical difference is that a compensation term is meaningful only within one uninterrupted
summation and becomes undefined when anything else modifies the state, whereas the pair represents
a value at every instant, which the reinitialization rule of Section~\ref{sec:full} relies on.

\subsection{Error bound}\label{sec:bound}

A bound follows from counting operations. Each two-float operation has a published worst-case
bound: barring overflow, one two-float addition (AccurateDWPlusDW) or multiplication (DWTimesDW) of
Joldes et al.\ returns a pair within $c\,\varepsilon^{2}$ of the exact result in relative terms, with
$c$ between $2$ and $6$ depending on the operation \cite{joldes2017,mullerrideau2022}. A Shu--Osher
row consists of a fixed number of these operations, and counting them for the rows of SSPRK(3,3)
gives an error below $30\varepsilon^{2}\approx10^{-13}$ per stage relative to
$\sum_k|\alpha_{ik}\hat u^{(k)}|+|\Delta t\,\beta_{ik}L_k|$, which for a variable with a large mean
is close to $|\bar u|$. The production product folds the low part of the leading product into the
first fused multiply--add; this variant saves one operation, but its bound is not among those proved
in \cite{joldes2017}. The split is exact and adds nothing. The two-float
state therefore introduces at most $30\varepsilon^{2}|\bar u|$ per stage, and at most about
$10^{-6}|\bar u|$ in total over $10^{7}$ stages, whereas an FP32 store can introduce
$\tfrac12\varepsilon|\bar u|$ at every stage, five orders of magnitude more: the factor of $2^{24}$
between the two unit roundoffs, less the operation count. Whatever the dynamics
do to an error once it has been introduced, they do equally to both representations. Apart from the
residual coefficient defect of $2^{-49}$, which is fixed in sign, these residuals do not have a
persistent sign, so in practice they accumulate
well below the worst case. In terms of Eq.~\eqref{eq:modeq}, the two-float state replaces the
forcing $R_u$ by a residual that is negligible compared with any practical truncation error, and
leaves $\tau$ and $R_L$ unchanged. Every operation is an FP32 addition, multiplication or fused
multiply--add, so that no FP64 instruction is executed and the cost does not depend on the FP64
rate of the hardware.

\subsection{Application to general schemes}\label{sec:storagerule}

A single rule covers every scheme of the form of Eq.~\eqref{eq:so}. Any vector that persists
across a stage boundary and holds values at the magnitude of the state is stored as a pair, whereas
registers that hold only increments stay in FP32; the rows are evaluated on $\hi+\lo$ and split on
store; $L$ always reads the $\hi$ float; and the values of $L$ are used within their own row, so
that they do not need a $\lo$ part. The memory cost depends on the sparsity of the coefficient
table, since $\lo$ floats are needed only for the state vectors that are live at the same time,
which are the registers a low-storage implementation already holds. For SSPRK(3,3) these are two
vectors, $u^{n}$ and the current stage, so the additional storage is two floats per equation, or
$40$ bytes per cell for a system of five equations. A low-storage scheme with three state
registers, such as the 3S* schemes of Ketcheson \cite{ketcheson2010}, requires $60$ additional bytes
per cell, and a Williamson 2N scheme \cite{williamson1980}, whose second register holds only
increments, requires $20$.

The rule does not require the scheme to be explicit, or to be a Runge--Kutta method. It requires
only that the update be an affine combination of stored states, with coefficients summing to one,
plus an increment. The increment form is the special case of a single coefficient equal to one. A
code that uses the Butcher form must therefore store as pairs the state vectors that it retains
across a step, which are fewer than in the Shu--Osher form. For BDF2, the retained time levels play
the role of the live vectors, with one $\lo$ float per stored level. A persistent per-step error is
amplified in BDF2 by the classical factor $1/\rho'(1)=3/2$, compared with at most one for the convex
rows of an SSP scheme, so BDF2 is at least as exposed to the state rounding error
(Fig.~\ref{fig:mech}d).

\subsection{Implementation in a production solver}\label{sec:full}

The two-float state was implemented in the Cartesian solver of the LAVA framework
\cite{nas_lava_cartesian}, a finite-volume solver for the compressible Navier--Stokes equations on
block-structured adaptively refined grids. In the original solver, the array kept between time
steps is the primitive state (pressure, velocity and temperature). At each stage the conservative
state is computed from it, updated, and converted back, so that the integrated quantity passes
through FP32 three times per step. In the modified solver the roles are reversed. The two-float
conservative state is kept between stages, because the Runge--Kutta rows combine conservative
variables, and the primitive state is recomputed from $\hi$ at every stage for the kernels that
use it.

One further convention keeps the change small. Parts of the code that modify the solution outside
the time integrator, such as mesh adaptation, load balancing and utilities, could be extended to
operate on the pair, but for convenience they instead mark the modification, so that the pair is
reinitialized once from the FP32 state and only a single store rounding is lost. Detecting such modifications by comparing the pair with the FP32 state was found
to be unreliable, because the conversion between primitive and conservative variables in FP32
introduces its own rounding errors of order one ulp, so that the comparison reports a modification
at nearly every stage and the $\lo$ floats are lost. Restart files contain both floats and have the
size of a double-precision restart file, and a continued run reproduces an uninterrupted one bit
for bit; FP32 restart files can still be read, in which case $\lo$ is initialized to zero. The
accessors use correctly rounded intrinsic functions and are therefore unaffected by the fast-math
compiler options that production GPU builds enable. When the two-float state is disabled, the code
compiles to the original solver, and in a double-precision build the accessors reduce to plain
loads and stores.

\section{Model problems}\label{sec:toys}

\begin{figure}[!tbp]\centering
\includegraphics[width=\textwidth]{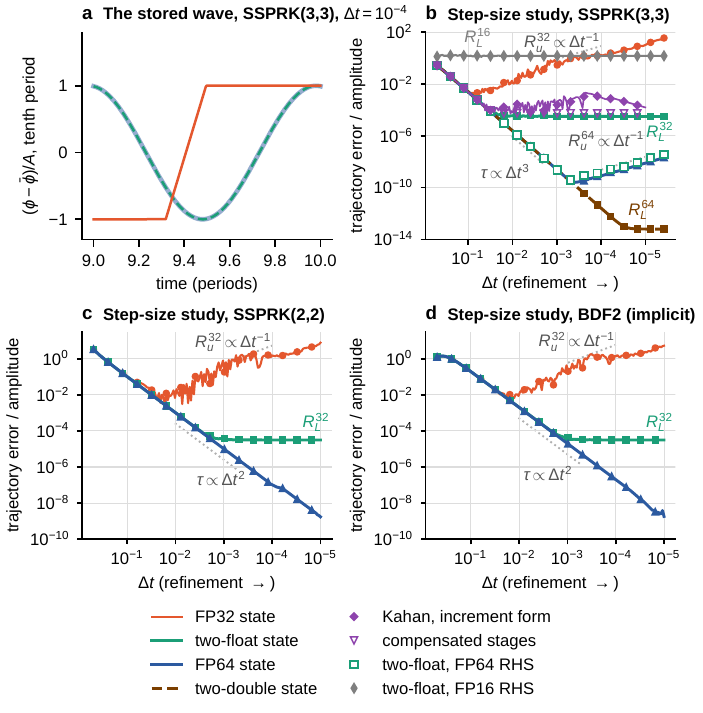}
\caption{Oscillator of amplitude $10^{-3}$ about a mean of $1.225$. (a) Stored solution over the tenth
period at $\Delta t=10^{-4}$ with SSPRK(3,3). (b) Trajectory error after ten periods as a function
of $\Delta t$ for SSPRK(3,3), at $171$ step sizes from $2^{-1}$ to $2^{-18}$, RMS over sixteen
initial phases. (c, d) As (b) for SSPRK(2,2) and implicit BDF2, at $157$ step sizes from $2^{-1}$ to
$10^{-5}$. Dotted lines mark the dominant term of Eq.~\eqref{eq:modeq} on each branch: the
truncation term $\tau$, the state rounding term $R_u\propto\Delta t^{-1}$ and the RHS rounding
floor $R_L$, with superscripts for the format of the state ($R^{32}_u$, $R^{64}_u$) or of the
right-hand side ($R^{32}_L$, $R^{16}_L$)}
\label{fig:mech}
\end{figure}

The model problems test the terms of the modified equation individually. Each problem is solved
with an FP32 state and the two-float state, which share all arithmetic and differ only in how the
state is stored between stages, and with an FP64 reference in which the state and all arithmetic,
including the right-hand side, are FP64; in the figures, ``FP64 state'' denotes this computation
entirely in FP64. The right-hand sides of the FP32 and two-float runs are evaluated in FP32 unless
stated otherwise. The two-float curve of Fig.~\ref{fig:mech}b is computed with FP32
arithmetic only. In the other panels of Fig.~\ref{fig:mech}, in its control runs, and in
Figs.~\ref{fig:order} and~\ref{fig:osc}, the pair is modelled by evaluating each row in FP64 on
$\hi+\lo$ and splitting the result, which agrees with the FP32-only curve to within $10\%$.

\subsection{Oscillator with a large mean}\label{sec:osc}

The simplest system that exhibits the behavior is a harmonic oscillation of small amplitude $A$
about a large mean $\bar\phi$,
\begin{equation}
\frac{d\phi}{dt}=\psi,\qquad \frac{d\psi}{dt}=-\omega^{2}(\phi-\bar\phi),
\label{eq:osc}
\end{equation}
with exact solution $\phi=\bar\phi+A\cos(\omega t+\theta)$. In every variant the restoring force is
evaluated in FP32 from the rounded value $\fl(\phi)$, in the same way that a flux kernel reads the
FP32 state.

\paragraph{Step-size study} Figure~\ref{fig:mech} shows the results of integrating
Eq.~\eqref{eq:osc} at the magnitudes of the density in a compressible flow, $\bar\phi=1.225$,
$A=10^{-3}$ and $\omega=1$, for ten periods. The error is measured in the $(\phi,\psi/\omega)$
phase plane and averaged over sixteen initial phases. At $\Delta t=10^{-4}$ the FP32 solution is
a staircase (Fig.~\ref{fig:mech}a): the stored value does not change on $82\%$ of the steps,
because near the turning points of the wave the stage increment is smaller than half an ulp of
the mean and is lost. As the step is refined (Fig.~\ref{fig:mech}b), the FP32 error reaches a
minimum of $1.7\times10^{-3}$ near $\Delta t=0.07$ and then increases as $\Delta t^{-1}$ for four
decades, until it exceeds the amplitude of the signal. This branch corresponds to $R_u$.

The two-float solution follows the FP64 solution along the truncation branch and levels off at
$3\times10^{-5}$, the floor $R_L$ set by the quantization of the FP32 state that the right-hand side
reads: the root-mean-square quantization error $\ulp(\bar\phi)/\sqrt{12}$ is
$3.4\times10^{-5}$ of the amplitude, and the measured floor is $3.1\times10^{-5}$ for every
integrator. In two control runs with the two-float state, an FP64
right-hand side makes it coincide with the FP64 solution, and an FP16 right-hand side raises its
floor to $1.5$ but leaves it flat: the precision of the right-hand side sets the level of the floor,
and the precision of the stored state decides whether the error grows.

FP64 is also affected. Its error turns around at $\Delta t=3.7\times10^{-4}$ and then increases as
$\Delta t^{-1}$, which is the same state rounding error at the FP64 ulp. For equal rounding
constants, the estimate of Section~\ref{sec:modeqn} predicts that the FP32 and FP64 turnaround steps
of this third-order scheme differ by the fourth root of the ratio of their unit roundoffs,
$(2^{29})^{1/4}\approx150$. The measured ratio is $0.072/3.7\times10^{-4}\approx190$; the FP32
minimum is shallow, and its location is uncertain by about a third. The same balance predicts a
ratio of minimum errors equal to the cube of the ratio of the steps, $7\times10^{6}$, compared with
a measured $6.6\times10^{6}$. Both measurements are consistent with an FP32 rounding constant about
$2.6$ times the FP64 one. The FP64 turnaround is pushed beneath the floor in the same way as the FP32
one: a state stored as two FP64 numbers with an FP64 right-hand side follows the $\Delta t^{3}$ branch a decade
beyond the FP64 turnaround and then levels off at $6.2\times10^{-14}$, five orders of magnitude
below the FP64 state at the smallest step. A wider number format thus
delays the turnaround, whereas widening only the stored state places the $\Delta t^{-1}$ term
beneath the right-hand-side floor over every step tested, at either precision. The pair has a
turnaround of its own, exposed by the control run with an FP64 right-hand side, which turns around
at the same step as the FP64 state (Fig.~\ref{fig:mech}b, open squares); above the floor, the
two-float state converges at the design order of the scheme. A Kahan-compensated update, the usual
single-float remedy, does not reach that floor here: once the increments fall below one ulp of the
stored density, the stage values formed from the FP32 state stagnate, and its error stays between
$10^{-4}$ and $10^{-3}$ (Section~\ref{sec:progression}). The second-order SSPRK(2,2) method
and the implicit BDF2 method, solved by fixed-point iteration to FP64 roundoff with the derivative
evaluated on the FP32 value of each iterate in every representation, show the same three regimes at their own
truncation orders (Fig.~\ref{fig:mech}c,d). For BDF2 the FP32 turnaround occurs at $\Delta t=2^{-6}$,
and the FP32 relative error reaches $1.2$ at $\Delta t=2^{-13}$, where the two-float state is at the
same floor as for the explicit schemes.

\begin{figure}[!tbp]\centering
\includegraphics[width=\textwidth]{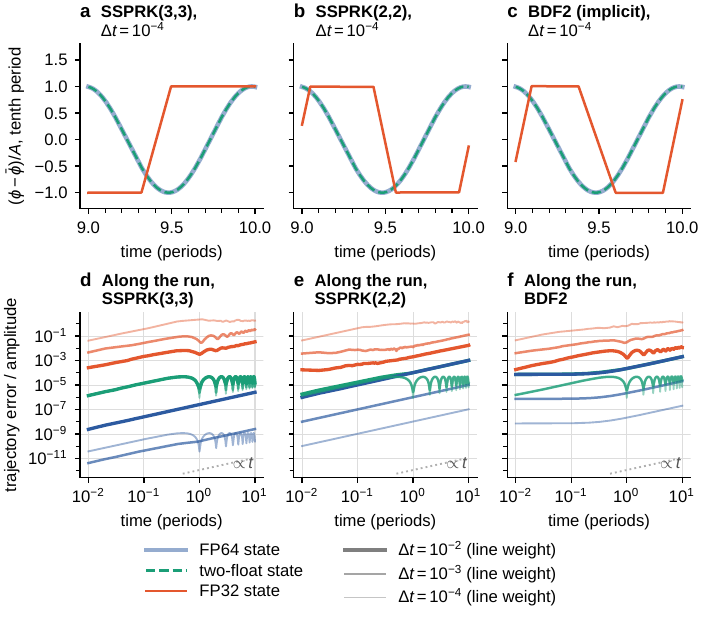}
\caption{Oscillator of Fig.~\ref{fig:mech} integrated with SSPRK(3,3), SSPRK(2,2) and implicit BDF2
(one column each). (a--c) Stored solution over the tenth period at $\Delta t=10^{-4}$. (d--f)
Trajectory error relative to the amplitude as a function of time at $\Delta t=10^{-2}$, $10^{-3}$
and $10^{-4}$, with heavier lines for larger steps (gray key below the panels), RMS over sixteen
initial phases, and a line of slope one for reference}
\label{fig:order}
\end{figure}

When the error is followed in time rather than evaluated at the final time (Fig.~\ref{fig:order}),
the FP32 solution is a staircase for every integrator, and its error grows in proportion to time
and is larger at smaller steps, while the two-float error remains at the floor set by the
right-hand side and the FP64 error, limited by truncation, decreases with refinement. The error is introduced
at every step from the first, so the problem is one of small time steps at a fixed horizon rather
than of long integration times.

\paragraph{A single step} The mechanism can be traced within one step of SSPRK(3,3) for
$dy/dt=r(y-Y_0)$ with $r\Delta t=0.01$ and $Y_0=10^{6}$, where the FP32 spacing is $0.0625$, starting
from $y^{0}=1{,}000{,}005.03$, which one float cannot represent ($\hi=1{,}000{,}005.0$,
$\lo=0.03$). The FP64 result is $y^{1}=1{,}000{,}005.0805523$. In FP32, the first stage is stored
with a rounding error of $+0.0125$, and the second-stage combination, $1{,}000{,}005.0283$, rounds
back to $1{,}000{,}005.0$, so that the whole stage increment is lost, as in the first experiment of
Section~\ref{sec:acc}; the error after one step is $-0.0181$, repeated at every step. The two-float
state uses the same FP32 stage derivatives but evaluates the combinations on $\hi+\lo$: its $\hi$
float moves to $1{,}000{,}005.0625$ at the first stage and stays there, while $\lo$ holds $+0.0175$,
$-0.0073$ and $+0.0180$ after the three stages. Its error after one step is $-3.15\times10^{-5}$,
$570$ times smaller and due entirely to the FP32 right-hand-side inputs, the $R_L$ term of
Eq.~\eqref{eq:modeq}, which does not grow under refinement.

\begin{figure}[!tbp]\centering
\includegraphics[width=\textwidth]{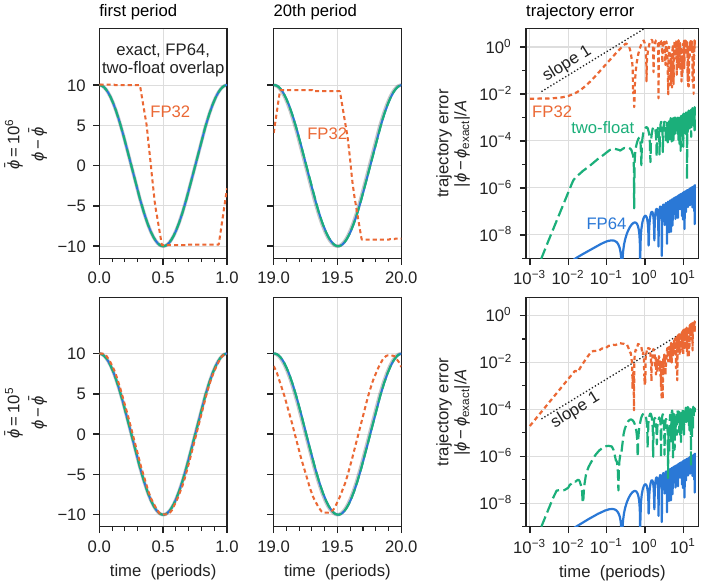}
\caption{The same oscillator about means of $10^{6}$ (top) and $10^{5}$ (bottom), integrated for
$60{,}300$ SSPRK(3,3) stages: the fluctuation over the first period (left) and the twentieth period
(middle), and the trajectory error relative to the amplitude over the whole run (right, logarithmic
axes, with a line of slope one for reference). On these axes an accumulating error appears as a
slope and a non-accumulating error as a horizontal line; the notches are zero crossings of the
error within a cycle. The FP64 and two-float solutions coincide with the exact solution at both
means}
\label{fig:osc}
\end{figure}

\paragraph{Dependence on the stored mean} To reach the regime of Section~\ref{sec:sign} in which
the increments are close to one ulp, the oscillator is placed at $\bar\phi=10^{6}$ with $A=10$,
$\omega=\ulp(\bar\phi)/A$ and $\Delta t=1$, and again one decade lower, where the same increments
span eight ulp. After twenty periods (Table~\ref{tab:osc} and Fig.~\ref{fig:osc}), the FP32 state
has lost $7.2\%$ and $2.3\%$ of its amplitude at the two means, and the two-float error is $200$ to
$500$ times smaller.

\begin{table}[!tbp]\centering\small
\caption{Final amplitude (exact value $10$) and its relative error, and the $L_1$ trajectory error
(the mean of $|\phi-\phi_{\mathrm{exact}}|$ over the run divided by $A$), after twenty periods of
the oscillator at a mean for which the increment per stage is one ulp, and at a mean one decade
lower. The amplitude is half the range of $\phi$ over the last period, so that distortion of the
waveform does not reduce the apparent error.}
\fitwidth{\begin{tabular}{lcccccc}\toprule
 & \multicolumn{3}{c}{$\bar\phi = 10^6$} & \multicolumn{3}{c}{$\bar\phi = 10^5$}\\
 \cmidrule(lr){2-4}\cmidrule(lr){5-7}
state & amplitude & error & $L_1$ & amplitude & error & $L_1$\\ \midrule
FP32 & $9.281$ & $-7.2\times10^{-2}$ & $0.90$ & $9.766$ & $-2.3\times10^{-2}$ & $0.16$\\
two-float & $10.00151$ & $+1.5\times10^{-4}$ & $1.1\times10^{-3}$ & $10.00119$ & $+1.2\times10^{-4}$ & $4.3\times10^{-5}$\\
FP64 & $9.99997$ & $-3.5\times10^{-6}$ & $4.1\times10^{-7}$ & $9.99997$ & $-3.5\times10^{-6}$ & $4.1\times10^{-7}$\\
\bottomrule
\end{tabular}}
\label{tab:osc}
\end{table}

\paragraph{Compressible state} Table~\ref{tab:fivevar} applies the scale rule to one oscillator
per conservative variable at sea-level, Mach $0.2$ magnitudes. In FP32 the variables with large
means, $\rho$ and $\rho E$, are in error by $13$ to $35\%$ and the momenta by less than one percent,
while the two-float errors are below $10^{-5}$. The rule identifies the exposed variables but does
not rank the others. The $\rho v$ row, which has no mean, suggests storing deviations from a
reference state (Section~\ref{sec:refstate}).

\begin{table}[!tbp]\centering\small
\caption{Scale rule applied to the compressible state: amplitude error, relative to the amplitude,
after $10^{5}$ SSPRK(3,3) steps for one oscillator per conservative variable at sea-level, Mach
$0.2$ magnitudes, with $\omega=1$ and $\Delta t=10^{-4}$ as in Fig.~\ref{fig:mech}; the error is that
of the invariant amplitude averaged over the final period ($\rho w$ behaves as $\rho v$ and is
omitted).}
\fitwidth{\begin{tabular}{lccccc}\toprule
variable & stored mean & fluctuation & $\ulp/\text{fluct.}$ & FP32 amp.\ error & two-float amp.\ error\\
\midrule
$\rho$    & $1.225$            & $10^{-3}$ & $1.2\times10^{-4}$ & $-0.35$ & $+9\times10^{-6}$\\
$\rho u$  & $83.4$             & $8.6$     & $8.9\times10^{-7}$ & $-1.6\times10^{-3}$ & $+7\times10^{-8}$\\
$\rho v$  & $\sim 0$           & $8.6$     & $1.1\times10^{-7}$ & $+2.4\times10^{-3}$ & $-2\times10^{-10}$\\
$\rho E$  & $2.561\times10^{5}$& $600$     & $2.6\times10^{-5}$ & $-0.13$ & $-1\times10^{-8}$\\
\bottomrule
\end{tabular}}
\label{tab:fivevar}
\end{table}

\subsection{Shallow water over a seamount}\label{sec:sw}

The shallow-water equations have the structure of an ocean or atmosphere model: a fast gravity
wave limits the time step, and a slow signal of a few centimeters is carried on a layer several
kilometers deep. We solve the one-dimensional nonlinear equations on a periodic domain of
length $L_x=2{,}000$\,km with $256$ cells ($\Delta x=7.8$\,km) over an ocean of depth
$H(x)=H_0-H_m\exp[-((x-L_x/2)/W)^{2}]$, with $H_0=4{,}000$\,m and a Gaussian seamount of height
$H_m=2{,}000$\,m and half-width $W=150$\,km. The stored variables are the free-surface elevation
$\eta_i=h_i-H_i$, where $h$ is the layer thickness, at cell centers and the velocity $u_{i+1/2}$ at faces (Arakawa C-grid), and the
semi-discrete equations are
\begin{equation}
\begin{aligned}
\frac{d\eta_i}{dt}&=-\frac{F_{i+1/2}-F_{i-1/2}}{\Delta x},\qquad
F_{i+1/2}=\tfrac12(h_i+h_{i+1})\,u_{i+1/2},\\
\frac{du_{i+1/2}}{dt}&=-g\,\frac{\eta_{i+1}-\eta_i}{\Delta x}
-u_{i+1/2}\,\frac{u_{i+3/2}-u_{i-1/2}}{2\Delta x},
\end{aligned}
\label{eq:sw}
\end{equation}
where $g$ is the gravitational acceleration and $h_i=\eta_i+H_i$ is computed in the working
precision at every evaluation, so that the
$4{,}000$\,m depth enters each right-hand side as an FP32 value but is never stored. The initial
condition is a right-going gravity wave of amplitude $a=0.1$\,m and wavelength $1{,}000$\,km,
integrated with SSPRK(3,3) in Shu--Osher form for $10^{4}$\,s. The reference solution is computed
in 80-bit arithmetic with $\Delta t=1.25$\,ms, and the error is the RMS of
$(\eta-\eta_{\mathrm{ref}})/a$ at the final time.

\begin{figure}[!tbp]\centering
\includegraphics[width=\textwidth]{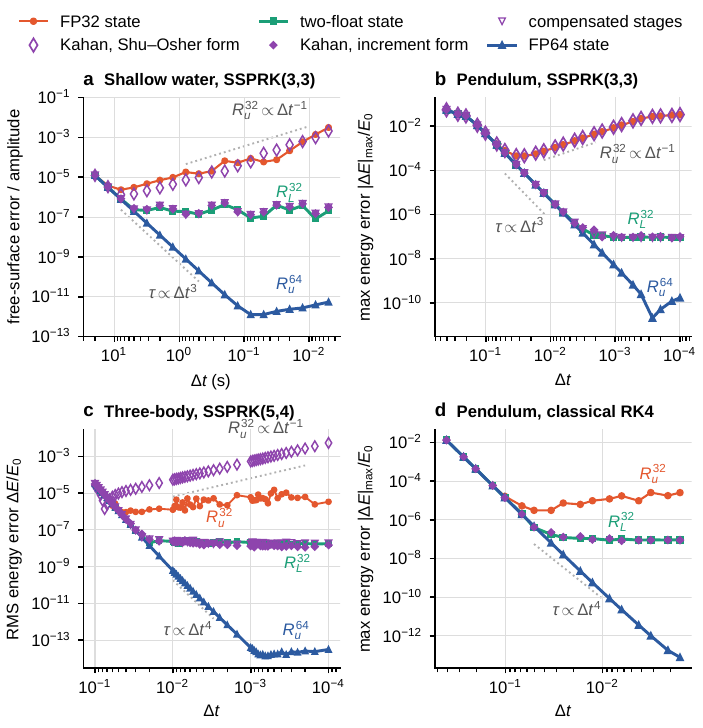}
\caption{Step-size studies at a fixed final time for three model problems. The FP32 and two-float
states share all arithmetic except the storage of the state between stages; the FP64 state is
computed entirely in FP64. (a) Nonlinear
shallow water over a seamount with SSPRK(3,3): RMS error of the free-surface elevation after
$10^{4}$\,s relative to an 80-bit reference. (b) Nonlinear pendulum released at $0.9\pi$ with
SSPRK(3,3): maximum energy error over $20$ periods, at $32$ to $131{,}072$ steps per period. (c)
Chaotic three-body figure-eight orbit with SSPRK(5,4) in increment form: RMS energy error over the
run to $T=10$, averaged over an ensemble of sixteen. (d) As (b) with classical RK4, whose
combination rows are all ones, so that the update is an increment: the Kahan-compensated update
and the two-float state coincide at every step size. Dotted lines mark the dominant term of
Eq.~\eqref{eq:modeq} on each branch, with symbols as in Fig.~\ref{fig:mech}}
\label{fig:domains}
\end{figure}

Storing the free-surface elevation rather than the layer thickness removes the large mean, so the
scale rule does not flag this problem, and the first coefficient of each row is computed as one
minus the rounded second, which is exact, so no coefficient defect enters. Nevertheless, the FP32
error turns around near the largest steps of the study, at
$\Delta t\approx8$\,s, and grows over three decades of refinement to $3\times10^{-3}$ of the wave
amplitude after two million steps, a factor of $1{,}300$ above its minimum
(Fig.~\ref{fig:domains}a). The fitted exponent is $-0.8$ over the first two decades and $-1.5$ over
the last $1.2$ decades, which is consistent with the intermediate case of Section~\ref{sec:sign}:
the fractional part of the increment drifts across the rounding grid by $a\omega^{2}\Delta t^{2}$
per step, so that a run lasts about $2\times10^{3}$ steps at $\Delta t=5$\,ms, and runs of that
length accumulate at a fixed horizon as $\Delta t^{-3/2}$. Kahan compensation applied to the same
Shu--Osher rows follows the FP32 error and reaches $2.1\times10^{-3}$ at the smallest step, because
the rounding of the products $\alpha_{ik}\hat u^{(k)}$ is lost before any addition it could protect
(Section~\ref{sec:progression}). Written in increment form, with a compensated update, the same
scheme instead follows the two-float state to its floor (filled diamonds).

A second program isolates the state rounding error by advancing each representation together with
an FP64 solution at the same $\Delta t$, so that the truncation errors cancel. At $\Delta t=5$\,ms
the FP32 error at the final time is $3.1\times10^{-3}$ of the amplitude and the two-float error is
$2.2\times10^{-7}$, smaller by a factor of fourteen thousand, although the two-float solver reads
the same FP32 state and differs only in how the state is stored. Storing only $\eta$ in FP32 gives
$2.1\times10^{-3}$, and storing only $u$ in FP32 gives $2.2\times10^{-3}$, so both variables
contribute. In time, the FP32 error first grows linearly and then approximately as the square root
of the number of steps, as expected for rounding errors that are correlated over a run. The
two-float solution follows the FP64 solution along the truncation branch and levels off at
$2\times10^{-7}$, which is the response to the FP32 rounding of the $4{,}000$\,m depth in the mass
flux. The FP64 error continues to decrease to $10^{-12}$ and then turns around as well, at $\Delta
t=0.08$\,s, although it grows more slowly than the FP32 error.

\subsection{Nonlinear pendulum and the scope of compensated summation}\label{sec:pendulum}

The nonlinear pendulum is a standard example of a geometrically nonlinear structural problem. With
$\ddot\theta+\sin\theta=0$ written as the system $u=(\theta,\omega)$, $L(u)=(\omega,-\sin\theta)$,
and released from rest at $\theta_0=0.9\pi$, the pendulum has $97.55\%$ of the energy $E=2$ needed
to reach the top, so an energy error of $2.45\%$ causes it to pass over the top. The reference is the exact
solution in terms of Jacobi elliptic functions, and the error is $\max|E-E_0|/E_0$, with
$E=\tfrac12\omega^{2}+2\sin^{2}(\theta/2)$ evaluated in FP64 from the stored state. One step of
SSPRK(3,3) is
\begin{equation}
\begin{aligned}
u^{(1)}&=u^{n}+\Delta t\,L(u^{n}),\\
u^{(2)}&=\tfrac34u^{n}+\tfrac14u^{(1)}+\tfrac14\Delta t\,L(u^{(1)}),\\
u^{n+1}&=\tfrac13u^{n}+\tfrac23u^{(2)}+\tfrac23\Delta t\,L(u^{(2)}),
\end{aligned}
\label{eq:pend}
\end{equation}
with each row evaluated and stored in the arithmetic of the representation under test; classical
RK4 and SSPRK(2,2) use the same interface.

Over $10{,}000$ periods at $1{,}024$ steps per period, the FP32 pendulum passes over the top at
period $765$, whereas the two-float pendulum completes the run and agrees with FP64 to three
digits. In a step-size study over twenty periods with up to $131{,}072$ steps per
period (Fig.~\ref{fig:domains}b), the FP32 energy error increases as $\Delta t^{-1}$ until the
pendulum passes over the top within the twenty periods and the error saturates at $3\%$. The
two-float error levels off at $10^{-7}$, the rounding level of the FP32 right-hand side, and the
FP64 error turns around at $\Delta t=3\times10^{-4}$.

The pendulum also provides the measurements behind the comparison with compensated summation in
Section~\ref{sec:progression}, which addresses a question raised by the geophysical studies cited in
Section~\ref{sec:prior} \cite{klower2022,lai2025}. For SSPRK(3,3) in Shu--Osher form, the
compensated pendulum passes over the top at period $723$, compared with $765$ for the FP32 state:
compensation changes which additions are rounded but not the outcome. With
classical RK4, whose steps are increments, compensation and the two-float state give the same
result at every step size (Fig.~\ref{fig:domains}d). Here the stage increments are many ulp, as they
are for the three-body problem of Section~\ref{sec:nbody}, so that the stage arguments formed from
the FP32 value without its compensation term do no harm. SSPRK(3,3) written in increment form
behaves in the same way: with a compensated update it follows the two-float state to its floor
(Fig.~\ref{fig:domains}b, filled diamonds). Applying the Shu--Osher coefficients in two-float form
does not rescue the Shu--Osher form: at $1{,}024$ steps per period it lowers the compensated error
about tenfold, to $6.5\times10^{-5}$, but the error still grows as $\Delta t^{-1}$, to
$5.9\times10^{-3}$ at the smallest step, because the rounding of each product
$\alpha_{ik}\hat u^{(k)}$ is still lost. Plain FP32 in increment form stalls between $10^{-5}$ and
$4\times10^{-5}$ instead of growing, which confirms that in this problem the coefficient defect
drives the growth of the Shu--Osher form. Over the $10{,}000$-period run,
the same scheme written in increment form completes the run with the FP64 energy error to within
$2\%$ for the FP32 state, the compensated update and compensated stages alike. In this problem,
where the update is the
entire computation, the two-float state requires $2.6$ times the FP32 time per step with RK4, and
compensation $1.2$ times.

\paragraph{Separating the two mechanisms}\label{sec:coeffs} Controlled experiments distinguish the two mechanisms of Section~\ref{sec:sign}. With increments
close to one ulp, as for the oscillator at $\bar\phi=10^{6}$, the state rounding dominates: a final
row of SSPRK(3,3) that sums to exactly one in FP32 lets the state stagnate completely, and the
amplitude error becomes $-100\%$ instead of $-7.2\%$. With increments of many ulp, as for the
pendulum, the coefficient defect dominates: a two-float state with
FP32 coefficients reproduces the FP32 error to within $3\%$, and only a two-float state and
two-float coefficients together gain four orders of magnitude. Classical RK4, whose rows are all
ones, has no coefficient defect.

\subsection{Chaotic three-body orbit}\label{sec:nbody}

The three-body problem provides a case in which the increments are many ulp. Three equal unit
masses move in the plane with $G=1$ on the figure-eight orbit of Chenciner and Montgomery
\cite{chenciner2000}, starting from $x_1=-x_2=(-0.97000436,\,0.24308753)$, $x_3=0$,
$v_1=v_2=(0.46620531,\,0.43236573)$, $v_3=-2v_1$, perturbed at the level of FP32 rounding to form an
ensemble of sixteen. The state is $u=(x_i,v_i)$ with $L(u)=(v_i,a_i)$,
$a_i=\sum_{j\ne i}(x_j-x_i)/|x_j-x_i|^{3}$, without softening. The five-stage fourth-order
SSPRK(5,4) method of Spiteri and Ruuth \cite{spiteri2002} is implemented in Butcher form,
\begin{equation}
y^{(i)}=u^{n}+\Delta t\sum_{j<i}a_{ij}\,L\big(y^{(j)}\big),\qquad
u^{n+1}=u^{n}+\Delta t\sum_{i}b_i\,L\big(y^{(i)}\big),
\label{eq:nbody}
\end{equation}
so that each step adds an increment to the state. The stages are evaluated in FP32 from the leading
float. The two-float state forms every stage argument and the update in two-float arithmetic, as in
Section~\ref{sec:rep}. For comparison, a Kahan-compensated increment applies TwoSum to the position
and velocity updates, which is sufficient for an increment-form update whose increments span many
ulp (Section~\ref{sec:progression}). The error is the energy error $(E-E_0)/E_0$, RMS over the steps of each run to $T=10$
and averaged over the ensemble.

The FP64 solution converges at order $4.26$ to its roundoff floor near $2\times10^{-14}$
(Fig.~\ref{fig:domains}c). The FP32 error stops decreasing between $10^{-6}$ and $10^{-5}$ after
the first decade of refinement, with a fitted exponent near zero. The two-float error levels off at about $2\times10^{-8}$ and
the compensated increment at about $1.4\times10^{-8}$, both near the level set by the FP32
right-hand side. The analysis explains why this exponent differs from the $\Delta t^{-1}$ of the
other model problems: with increments of many ulp, the state rounding errors do not keep a
consistent sign, and an increment-form update has no coefficient defect, so the error stagnates
rather than growing. Written in Shu--Osher form, whose FP32 combination rows sum to one
only within $2^{-25}$, the same scheme behaves like the other model problems: the FP32 error grows
as $\Delta t^{-1}$, to $1.8\times10^{-3}$ at $\Delta t=3\times10^{-4}$, and Kahan compensation
follows it (Fig.~\ref{fig:domains}c, open diamonds). In neither form does the FP32 solution
converge.

\subsection{Alternative remedies}\label{sec:refstate}

Since the bias originates in the store of each stage, a less expensive remedy than a second float
might suffice, and we examined the remedies that are commonly proposed. Nondimensionalization,
storage relative to a reference state and the promotion of selected arithmetic to FP64 change the
solver outside the time integrator, in its variables, boundary conditions, output or flux kernels.
The others change only the time integrator, so that output, mesh adaptation and restart are
unaffected; Table~\ref{tab:remedies} compares them for four integrators, each in its usual form.

Nondimensionalization makes the stored means of order one so that their ulps appear small, but
floating-point numbers have relative precision and the ratio of Eq.~\eqref{eq:scale} is
dimensionless. Rescaling the FP32 oscillator by $2^{20}$
reproduces its trajectory bit for bit, and the density, which is stored at order one, is the most
affected variable in Table~\ref{tab:fivevar}. A factor that is not a power of two moves the stored
mean between consecutive powers of two and can make the error larger: in the wing-station problem
below, scaling the pressure by $\rho_\infty a_\infty^{2}$ raises the error at $\Delta t/4$ from
$+127$--$147\%$ to $+552$--$561\%$. Reducing the stored magnitude requires subtracting
an offset rather than rescaling.

Storing the deviation from a uniform reference state, for example a gauge pressure
$p-p_{\mathrm{ref}}$ and a temperature difference, has been used in
reduced-precision lattice Boltzmann codes \cite{lehmann2022}, and a related approach, in which the
reference matches the local mean everywhere, is the classical remedy for cancellation at low Mach
number \cite{sesterhenn1999}. With a uniform reference, the stored value is rounded at
$\ulp(|\Delta|+A)$, where $\Delta$ is the offset of the local mean from the reference, and the flow
places the largest offsets at the stations that generate the loads. The model problem consists of
five copies of the oscillator at the pressure offsets of five wing stations, $C_p\,q_\infty$ with
$C_p\in\{0,+1,-1,-4,-10\}$ (far field, stagnation, mild suction, suction peak and slat suction),
with identical fluctuations ($A=10$\,Pa), an increment per stage of one ulp of the absolute
pressure, and the most favorable implementation of the reference-state scheme, in which all
arithmetic is carried out on the deviation and the absolute pressure is never formed.
Table~\ref{tab:stations} gives the amplitude error after twenty periods at $\Delta t$ and ten periods
at $\Delta t/4$. At the loaded stations the reference state makes the ulp of the stored value only
$4$ to $32$ times smaller, which delays the onset of the error by the same factor in $\Delta t$ but
does not remove it; refinement brings it back, first at the suction peaks. Reference-state storage
is nevertheless useful together with the two-float state, for the reason given in
Section~\ref{sec:guidance}.

\begin{table}[!tbp]\centering\small
\caption{Reference-state storage with identical fluctuations at five wing stations, after twenty periods
at $\Delta t$ and ten periods at $\Delta t/4$. For comparison,
absolute FP32 storage gives errors of $+12$ to $+41\%$ at $\Delta t$ and $+127$ to $+147\%$ at
$\Delta t/4$ at every station, and the two-float state gives $0.02\%$ or less everywhere. The error
of the reference-state scheme follows the increment measured in ulp of the stored deviation.}
\fitwidth{\begin{tabular}{lrrrrr}\toprule
& & \multicolumn{2}{c}{increment / $\ulp$(stored deviation)} & \multicolumn{2}{c}{reference-state error}\\
\cmidrule(lr){3-4}\cmidrule(lr){5-6}
station & $\Delta$ (Pa) & at $\Delta t$ & at $\Delta t/4$ & at $\Delta t$ & at $\Delta t/4$\\
\midrule
far field ($C_p=0$)        & $0$       & $8192$ & $2048$ & $+0.5\%$ & $+0.9\%$\\
stagnation ($C_p=+1$)      & $+2837$   & $32$   & $8$    & $-0.2\%$ & $+1.7\%$\\
mild suction ($C_p=-1$)    & $-2837$   & $32$   & $8$    & $-0.2\%$ & $+1.7\%$\\
suction peak ($C_p=-4$)    & $-11348$  & $8$    & $2$    & $+1.3\%$ & $+16.3\%$\\
slat suction ($C_p=-10$)   & $-28371$  & $4$    & $1$    & $+4.1\%$ & $+54.2\%$\\
\bottomrule
\end{tabular}}
\label{tab:stations}
\end{table}

Promoting selected arithmetic to FP64 while storing the state in FP32 leaves the store rounding
unchanged; four such promotions in the production solver gave no improvement in accuracy and
increased the cost by $8.4\%$.

Table~\ref{tab:remedies} compares the remedies that change only the time integrator on the
wing-station problem at $\Delta t/4$, with SSPRK(3,3), classical RK4, the ten-stage, fourth-order
SSPRK(10,4) in the two-register low-storage form of Ketcheson \cite{ketcheson2008}, and BDF2, each in
its usual form, since rewriting an integrator in another form is not a drop-in change. With SSPRK(10,4), whose low-storage form
multiplies one register by $15$ and another by $5$ before they cancel, FP32 reaches $+1757\%$.
Compensation applied in the usual form fails wherever coefficients multiply the state; it helps only
for RK4, whose usual form is already the increment form, and compensated stages then reach the
floor. The two-float state reaches the floor in the usual form of every integrator. Rewriting the
integrator in increment form, which is not a drop-in change for SSPRK(3,3), SSPRK(10,4) or BDF2,
alters this picture: compensated stages then reach the floor with every integrator, at the cost of
the rewrite and, for SSPRK(10,4), of storing up to ten stage derivatives in place of two registers,
whereas a delta-form update in FP32 lowers the error only to $+51\%$, and the FP32 difference form of
BDF2 raises it to $+3707\%$. Stochastic
rounding \cite{croci2023} removes the bias of the store rounding but not its variance, and does not
retain the rounding error; it is a natural alternative where the hardware supports it, and a
quantitative comparison is left for future work.

\begin{table}[!tbp]\centering\small
\caption{Remedies that change only the time integrator, for four integrators: amplitude error after
ten periods at $\Delta t/4$ for the five wing-station oscillators of Table~\ref{tab:stations}, as
the range over the stations in percent of the FP64 amplitude of the same integrator. Each integrator
is used in its usual form (Shu--Osher for SSPRK(3,3), Butcher for RK4, the two-register low-storage
form \cite{ketcheson2008} for SSPRK(10,4), and $u^{n+1}=\tfrac43u^{n}-\tfrac13u^{n-1}+\ldots$ for
BDF2). Every state is read at its full value: $\hi+\lo$ for the two-float state and the value with
its correction for the compensated variants. A dash marks a remedy that requires rewriting the
integrator in another form; for RK4, whose usual form is the increment form, Kahan compensation is
the compensated update, and all coefficients that multiply the state are one (n/a).}
\footnotesize\setlength{\tabcolsep}{4pt}
\fitwidth{\begin{tabular}{lrrrr}\toprule
remedy & SSPRK(3,3) & RK4 & SSPRK(10,4) & BDF2\\ \midrule
FP32 state & $+127$ to $+147$ & $+51$ & $+205$ to $+1757$ & $+93$ to $+542$\\
Kahan, usual form & $+47$ to $+369$ & $+0.30$ & $+413$ to $+1273$ & $+478$ to $+1063$\\
Kahan with two-float $\alpha$, usual form & $+77$ to $+282$ & n/a & $+86$ to $+821$ & $+468$ to $+1096$\\
compensated stages & --- & $-0.001$ & --- & ---\\
two-float state, usual form & $-0.02$ to $+0.01$ & $-0.02$ to $+0.01$ & $-0.02$ to $+0.02$ & $-0.02$ to $+0.01$\\
\bottomrule
\end{tabular}}
\label{tab:remedies}
\end{table}

\section{Application to a production flow solver}\label{sec:production}

\subsection{Solver and test cases}\label{sec:setup}

The production solver integrates the compressible Navier--Stokes equations in the conservative
form
\begin{equation}
\partial_t W + \nabla\!\cdot\!F_c(W) = \nabla\!\cdot\!F_v(W,\nabla W),
\qquad W=(\rho,\;\rho u,\;\rho v,\;\rho w,\;\rho E)^{T},
\label{eq:ns}
\end{equation}
where $F_c$ and $F_v$ are the convective and viscous fluxes, with the state $W$ stored in SI
units, on block-structured adaptively refined Cartesian grids \cite{nas_lava_cartesian}. The
solver uses SSPRK(3,3) at a CFL number near one, a central scheme with shock and turbulence
sensors, the Vreman subgrid-scale model \cite{vreman2004}, an immersed boundary method, and an
equilibrium wall-stress model with a transition treatment, and runs on GH200 GPU nodes \cite{nvidia2025gh200}. Three
builds from the same version of the solver are compared: FP32, FP64, and the two-float state of
Section~\ref{sec:method}. The FP32 and FP64 builds are the original solver, which keeps the primitive
state between stages in its own precision and converts to and from the conservative variables at
every stage (Section~\ref{sec:full}); the two-float build keeps the conservative state as pairs. The
builds therefore differ in the stored variables as well as in their precision. The FP64 and
two-float builds, which store different variables, agree on the loads, whereas the FP32 build, which
stores the same variables as the FP64 build, does not (Section~\ref{sec:hlpw}), so the choice of
stored variables does not explain the FP32 bias. All three are compiled with the same
options, including the fast-math options that the solver's GPU builds enable, and compared runs
use identical input files.

\subsection{Taylor--Green vortex}\label{sec:tgv}

The first test case is the inviscid Taylor--Green vortex \cite{taylor1937} at a Mach number of
about $0.1$ on a periodic $32^{3}$ grid with a single level, run for ten convective times ($5{,}120$
steps at the base $\Delta t$). The inviscid equations conserve entropy, so the relative increase of
the integrated entropy over the run measures the numerical dissipation in a single number. For the
FP64 build this increase is $1.179\times10^{-4}$, unchanged to four digits when $\Delta t$ is
halved and quartered, so it is due to the spatial truncation error. The FP32 build gives
$6.40\times10^{-4}$, $1.17\times10^{-3}$ and $2.21\times10^{-3}$ at $\Delta t$, $\Delta t/2$ and
$\Delta t/4$. Subtracting the FP64 value leaves an FP32 excess of $5.22\times10^{-4}$,
$1.05\times10^{-3}$ and $2.09\times10^{-3}$, which grows by factors of $2.01$ and $1.99$ each time
the step is halved; the $\Delta t^{-1}$ dependence of Eq.~\eqref{eq:modeq}
therefore holds to within one percent in the full solver. The difference in kinetic energy shows the
same trend less clearly, $7.7\times10^{-4}$, $1.11\times10^{-3}$ and $2.50\times10^{-3}$, because it
combines the effect of the store rounding with other differences between the builds. The two-float
build gives $1.178\times10^{-4}$ at the base step and $1.182\times10^{-4}$ at $\Delta t/4$, equal to
the truncation contribution, while the FP32 excess grows to eighteen times that value. An FP32 build
without fast-math options reproduces the FP32 results to within a few percent ($6.50\times10^{-4}$
and $2.21\times10^{-3}$). For the pressure field, the two-float state reduces the error relative to
FP64 by a factor of $427$ at $512$ steps and keeps it between $0.24$ and $0.57\%$ up to $20{,}480$
steps, whereas the FP32 error continues to grow.

\subsection{Wall-modeled flow past a sphere}\label{sec:sphere}

The second test case exercises all production components at once, including the immersed boundary
that represents the body on the Cartesian mesh, the wall model, and the surface integration of the
loads, and measures the error introduced per step. All three builds start from the same FP32 field
of a developed wall-modeled flow past a sphere; since conversion from FP32 to FP64 is exact, the
initial states are identical. Each build advances $K=1$, $10$ and $100$ steps at a fixed $\Delta t$,
and the volume solution is compared with the FP64 build for each component in the maximum norm,
normalized by the range of the component. After one step the FP32 build reproduces the ranking of
Table~\ref{tab:fivevar}, with errors of $2.9\times10^{-5}$ in temperature, $9.2\times10^{-6}$ in
pressure and $5$--$9\times10^{-7}$ in the velocity components. By $K=100$ the FP32 error has grown by
a factor of nine to twenty-two, depending on the component, while the two-float error remains within
a factor of two of its value after one step. These factors are ratios of total errors, and the
one-step error already contains RHS rounding errors that do not grow, so the growth of the state
rounding error alone is larger than the ratios suggest. In addition, the maximum norm over a
heterogeneous field is not the quantity governed by the $\Delta t^{-1}$ dependence, which applies to
the phase-correlated mean that the loads integrate and that Section~\ref{sec:tgv} measures. The
spatial distribution of the error is consistent with the analysis: the largest FP32 errors in
pressure and temperature occur two to three diameters from the body, in the undisturbed stream,
where the stored state consists almost entirely of its mean.

\subsection{High-lift wing}\label{sec:hlpw}

The main test case is Test Case~1 of the Sixth AIAA High-Lift Prediction Workshop, the CRM-HLS
\cite{hlpw6cases}, a swept, finite-span, constant-section high-lift wing with slat, flap and slat
brackets, designed to reproduce the outboard separation of the full high-lift Common Research Model
at reduced cost. At a Mach number of $0.2$, a Reynolds number of $3.55\times10^{6}$ and an angle of
attack of $30^\circ$, the outboard separation that limits the maximum lift is present. The case is
sensitive to the state rounding error for three reasons: the stored means are large compared with the fluctuations of interest, an averaging
window contains millions of stage updates, and the quantities of interest are surface-integrated
coefficients reported to one part in $10^{4}$.

The three builds are started from the same initial condition and run on the same decomposition over
sixteen GPUs throughout the refinement schedule, so that differences in decomposition cannot be
confused with differences in precision. Grid levels are added at $50$, $100$ and $150$ convective
time units (CTU), each halving the finest cell size and, at fixed CFL number, the time step. The
finest cell size at level $\ell$ is $\Delta x=L_{\mathrm{ref}}/2^{\ell}$ with
$L_{\mathrm{ref}}=0.762$\,m, so the three refined grids, at levels $9$, $10$ and $11$, have finest
cells of $1.5$, $0.74$ and $0.37$\,mm and contain $60$, $230$ and $544$ million cells. The number of
steps per CTU roughly doubles at each refinement, to about $10{,}400$, $23{,}500$ and $47{,}000$. Because the refinement is
solution-adaptive, the number of cells grows by less than a factor of eight per level. Averaging
windows begin ten CTU after the first two refinements and five CTU after the last, and end at the
next refinement or at the end of the shortest run, giving windows of $60$--$100$, $110$--$150$ and
$155$--$171$ CTU. The FP64 build serves as the reference; the comparison therefore measures
consistency with FP64, not agreement with experiment.

Figure~\ref{fig:trio}, shown in the introduction, gives the coefficient histories, and Table~\ref{tab:loads} gives the
statistics. The FP32 build overpredicts the lift by $2.1\%$, $4.5\%$ and $7.7\%$ of the FP64 value
on the three grids, and the drag by $1.6\%$, $4.1\%$ and $6.9\%$. The difference from FP64 roughly doubles with each refinement, as does the
number of steps per convective time: refining the grid, which normally brings a solution closer to
the converged answer, moves the FP32 solution further from the FP64 solution on the same grid.
Because the grid and the time step are refined together, this comparison does not separate their
effects; the Taylor--Green test of Section~\ref{sec:tgv}, which halves $\Delta t$ on a fixed grid,
isolates the dependence on the time step. Within the two finer windows the FP32 lift is still rising
(Fig.~\ref{fig:trio}), so for those windows its percentages may understate the eventual difference.
The two-float results agree with the reference: the lift differences are
$+0.26$, $0.00$ and $-0.19\%$, and on the middle window, where FP32 is $4.5\%$ high in lift and
$4.1\%$ high in drag, the two-float differences are $0.00\%$ and $-0.01\%$.

\begin{table}[!tbp]\centering\small
\caption{Differences in the load coefficients at $30^\circ$ angle of attack, in counts ($10^{-4}$)
and as a percentage of the FP64 window mean, with $95\%$ confidence intervals, for one
campaign started from a single initial condition and run on one sixteen-GPU decomposition through
all refinement stages.}
\footnotesize\setlength{\tabcolsep}{3pt}
\fitwidth{\begin{tabular}{llcrrrrr}\toprule
 & grid & window (CTU) & FP64 mean & FP32$-$FP64 & (\%) & two-float$-$FP64 & (\%)\\ \midrule
$\cl$  & $60$M  & $60$--$100$  & $1.9085$ & $+409\pm193$   & $+2.1\pm1.0$ & $+49\pm148$    & $+0.26\pm0.77$\\
       & $230$M & $110$--$150$ & $2.2257$ & $+992\pm208$   & $+4.5\pm0.9$ & $+0.6\pm74.0$  & $+0.00\pm0.33$\\
       & $544$M & $155$--$171$ & $2.3742$ & $+1835\pm413$  & $+7.7\pm1.7$ & $-45\pm69$     & $-0.19\pm0.29$\\
\midrule
$\cd$  & $60$M  & $60$--$100$  & $0.4149$ & $+66\pm60$     & $+1.6\pm1.5$ & $-22\pm46$     & $-0.54\pm1.11$\\
       & $230$M & $110$--$150$ & $0.4440$ & $+182\pm16$    & $+4.1\pm0.4$ & $-0.6\pm7.8$   & $-0.01\pm0.18$\\
       & $544$M & $155$--$171$ & $0.4697$ & $+323\pm64$    & $+6.9\pm1.4$ & $-12.8\pm7.7$  & $-0.27\pm0.16$\\
\bottomrule
\end{tabular}}
\label{tab:loads}
\end{table}

The two-float results also bound the contribution of the solution-adaptive refinement, through
which the state rounding error could reach the loads by way of the mesh. The two-float solution
departs from the FP64 solution at the level of its own rounding residual from the first step, so
its refinements, like those of the FP32 build, respond to a solution decorrelated from the
reference; if differences in the grid alone could change a window mean by several percent, they
would do so in the two-float run as well, and they do not. The same argument answers the objection that any two decorrelated turbulent runs
might differ by this amount: the two-float and FP64 runs are two such runs, and they agree.

\begin{figure}[!tbp]\centering
\includegraphics[width=\textwidth]{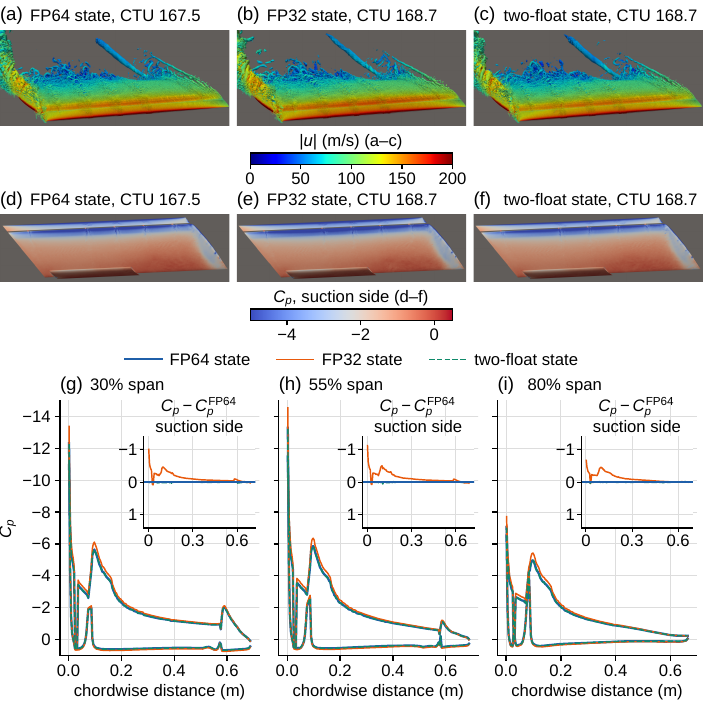}
\caption{HLPW-6 Test Case~1 (CRM-HLS) at $30^\circ$ angle of attack on the $544$M-cell grid near the
end of the campaign, one column per build: FP64 at $167.5$ CTU, and FP32 and two-float at $168.7$
CTU. (a--c) Isosurfaces of the $Q$-criterion at $10^{5}$\,s$^{-2}$ colored by velocity magnitude.
(d--f) Instantaneous pressure coefficient on the suction side, viewed normal to the planform. (g--i)
Time-averaged surface pressure (the solver's five-CTU exponential moving average) at three spanwise
stations, with the builds overlaid and $C_p$ plotted with suction upward; the insets show the
suction-side difference from FP64 on a scale ten times finer. Each row uses one viewpoint and one
color range}
\label{fig:edtrio}
\end{figure}

The bias is not visible in instantaneous fields, but it is clear in the statistics.
Figure~\ref{fig:edtrio} shows the three builds near the end of the campaign. Their isosurfaces and
instantaneous surface pressures are decorrelated realizations of the same turbulent flow and cannot
be distinguished. The time-averaged surface pressure, however, separates them. The FP32 suction
peaks exceed those of FP64 by $0.66$ to $1.11$ in $C_p$ at the slat and by $0.46$ to $0.49$ on the
main element, and the mean suction-side $C_p$ of FP32 is $0.12$ to $0.13$ more negative at every
station, which is the surface signature of the lift difference. The two-float results agree with
FP64 to within $0.03$ at the peaks and $0.007$ in the mean.

\paragraph{Other possible sources of the bias}\label{sec:exon} The model problems identify the store of the state as the source of the bias, and the production
solver confirms this by elimination. Each other possible source was tested with a separate
instrumented build, run against the standard FP32 build in short matched runs from a common state
and compared step by step, so that chaotic divergence does not mask small differences. Removing the fast-math
intrinsic functions and promoting the flux and right-hand-side assembly to FP64 left the result
unchanged. Disabling the wall model did not change the divergence between the builds. Running the
subgrid-scale model and the shock and turbulence sensors in FP64 changed the window-mean loads by
less than $0.01$ counts, and computing the entire load integration in FP64 reproduced the standard
FP32 result to better than $0.1$ count. Within the accuracy of these measurements, the store of the
state is the only biased term we have found.

\subsection{Cost}\label{sec:cost}

Table~\ref{tab:cost} gives the measured cost, normalized by that of the FP32 build, for three
configurations: an inviscid Taylor--Green vortex on $256^{3}$ cells and a WMLES of the HLPW-5
high-lift model at $19.7^\circ$ angle of attack on $14.9$ million cells, each on one GH200 GPU, on
which FP64 runs at half the FP32 rate, and the same Taylor--Green case on an RTX Blackwell
workstation GPU, on which FP64 runs at one sixty-fourth of the FP32 rate. On the GH200, the
two-float build takes $16\%$ more time per step than the FP32 build for the Taylor--Green vortex
and $12\%$ more for the WMLES, whereas the FP64 build takes $3.1$ and $1.8$ times as long. The
two-float build is therefore $2.7$ and $1.6$ times faster than FP64 for these two cases. For the
WMLES it also needs $7\%$ more memory than FP32, compared with $58\%$ more for FP64. The FP64
memory ratio is below two because grids, masks and connectivity data do not depend
on the precision. Of the additional $7\%$ of the two-float build, the $40$ bytes per cell of
Section~\ref{sec:storagerule} account for about four fifths, and halo cells and temporary arrays
for the remainder. Even with FP64 at half rate, most of the cost of
double precision is due to the doubled memory traffic and footprint, which faster FP64 arithmetic
does not reduce.

On the workstation GPU the cost of the two-float build is unchanged at $1.17$, whereas double
precision is not a practical option on that hardware. The cost does not depend on the FP64 rate of
the hardware because the accessors execute no FP64 operations. An earlier implementation that used
FP64 arithmetic inside the accessors cost $1.26$ on the same GPU. That implementation corresponds to
the mixed-precision design of OpenMM \cite{eastman2017}, which stores the same two floats but updates
them in FP64 arithmetic.

Another alternative is an FP64 state with every right-hand-side evaluation left in FP32, as in
the mixed-precision schemes of Section~\ref{sec:prior}. It can be given the same structure as the
two-float state: the FP64 state is updated at each store, and its rounded value is written into the
existing FP32 array that every other kernel reads, so the flux kernels are again unchanged. We did
not implement this variant, but its store performs FP64 arithmetic as the FP64-promoted accessors
do, so on the workstation GPU it should cost at least the $1.26$ measured for them, against $1.17$
for the two-float state, and it holds twelve bytes per stored variable, the FP64 state and its FP32
image, against eight for the pair. Where FP64 arithmetic is fast, the two designs should cost about
the same apart from that memory; where it is slow, only the two-float state keeps its cost.
Because the two-float arithmetic is confined to the store, a few dozen FP32 operations per cell and
stage, and the flux kernels are unchanged, the earlier conclusion that emulated precision is too
expensive \cite{goddeke2007} does not apply here. The remaining cost is the memory traffic and
footprint of the $\lo$ floats.

\begin{table}[!tbp]\centering\small
\caption{Measured cost, normalized by the FP32 build (larger values are worse). The timing rows are
the reciprocal of the throughput. The FP64-promoted row uses the same two-float scheme with FP64
arithmetic inside the accessors.}
\fitwidth{\begin{tabular}{lccc}\toprule
measurement & FP32 & two-float & FP64\\ \midrule
time per step, TGV $256^3$, 1 GPU (GH200, FP64 $1{:}2$) & $1.00$ & $1.16$ & $3.13$\\
time per step, HLPW-5 WMLES $14.9$M cells, 1 GPU (GH200) & $1.00$ & $1.12$ & $1.82$\\
time per step, TGV $256^3$, 1 GPU (RTX Blackwell, FP64 $1{:}64$) & $1.00$ & $1.17$ & ---\\
\quad the same, with FP64-promoted accessors & --- & $1.26$ & ---\\
GPU memory, HLPW-5 WMLES $14.9$M cells & $1.00$ & $1.07$ & $1.58$\\
\bottomrule
\end{tabular}}
\label{tab:cost}
\end{table}

\section{Discussion}\label{sec:discussion}

\subsection{Relation to earlier studies}\label{sec:earlier}

The analysis of Section~\ref{sec:modeq} reconciles the present results with the earlier findings
that single precision is adequate \cite{homann2007,witherden2020,vana2017,karp2025}. The most direct
comparison is with Karp et al.\ \cite{karp2025}, who rounded the stored state as well as individual
terms. With the state
rounded to FP32 their solvers agreed with FP64, whereas with the state rounded to sixteen bits or
fewer several runs showed the behavior described in Section~\ref{sec:acc}, including a channel flow
at a friction Reynolds number of $683$ instead of $550$. Three differences account for why FP32 was
sufficient in those studies but not here. The first applies to the incompressible solvers; the second
and third apply to the compressible ones as well.

The first is the stored quantity. In incompressible flow the mean velocity is at most the bulk
velocity and the turbulent fluctuations are a substantial fraction of it, so the ratio of
Eq.~\eqref{eq:scale} is near $10^{-6}$ in FP32, $25$ to $120$ times smaller than for the energy and
density of a compressible flow at Mach $0.2$ (Table~\ref{tab:fivevar}). The second is the time
integrator. The classical fourth-order Runge--Kutta method, used in Simson, one of their incompressible solvers, has combination rows
that are all ones and therefore no coefficient defect, which Section~\ref{sec:coeffs} identified as
the systematic error when the increments are many ulp, and the Bogacki--Shampine scheme of SSDC, their compressible solver, in
Butcher form, carries fewer roundings with weight one to the end of a step than a scheme in
Shu--Osher form. The third is the quantity of interest. Their cylinder statistics were averaged over
about a hundred convective times, and since their FP64 solvers differ from one another by several
percent on these statistics, differences of the size in Table~\ref{tab:loads} would not be visible.
A load coefficient reported to $10^{-4}$ is more demanding, and one averaging window of the wing
campaign contains more than two million stage updates per variable. For the compressible
flux-reconstruction simulations of Witherden and Jameson \cite{witherden2020}, the closest
counterpart to the present ones, the second and third differences likely apply: their low-storage
scheme updates the state in increment form, a high-order discretization takes far fewer steps per
flow-through than a second-order one, and they compared turbulence statistics rather than loads
reported to one part in $10^{4}$.

The analysis also bears on how precision should be chosen. It has been proposed that precision be
treated as one more simulation parameter, to be balanced against resolution and time step
\cite{karp2025}. For the stored state the two are not independent, because the error is driven by a
forcing of magnitude $\ulp(\bar u)/\Delta t$: refining the time step increases it, and a precision
found adequate on a coarse grid can become inadequate on a refined one. Accumulators such as the
simulation time and the weights of running statistics, which are commonly kept in FP64 for the same
reason \cite{karp2025}, are of the same kind as the stored state, and the two-float state
applies the same remedy to the state without requiring FP64.

\subsection{Recommendations}\label{sec:guidance}

Whether a code is affected can be determined at low cost. The time step is halved for a case with a
fixed final time, and the result is compared with a double-precision run. For a chaotic flow the
comparison must use a quantity that does not depend on the trajectory, such as an integrated
dissipation or a time-averaged load. If the difference grows,
and the simulation time and statistics accumulators are already held in higher precision, the store
of the state is the cause. The Taylor--Green results of Section~\ref{sec:tgv}, with growth factors of
$2.01$ and $1.99$, illustrate the outcome for an affected code. The scale rule of Eq.~\eqref{eq:scale}, evaluated
from the stored magnitudes and the signal of interest without reference to the discretization,
indicates which variables are affected and how the exposure changes with the choice of units,
reference state or resolution.

The remedy does not depend on the physics of the problem. The storage rule of
Section~\ref{sec:storagerule} applies, the combination coefficients are held in two-float form, and,
unless it is extended to the pair, any modification of the solution outside the time integrator is
marked so that the pair is reinitialized once rather than discarded (Section~\ref{sec:full}). For an
integrator written in increment form the second number can be a compensation term, provided that it
also enters every stage value; where scalars multiply the state, as in the Shu--Osher and multistep
forms, it must be a value that is closed under multiplication. The memory cost is one additional
number per equation for each live state vector, two for SSPRK(3,3) and for the two-register
SSPRK(10,4), and the cost in time is $12$--$17\%$ for a solver dominated by RHS evaluation. Where
the update is the entire computation, as for the pendulum, the two-float state costs $2.6$ times as
much as FP32; such problems are usually inexpensive, and where FP64 arithmetic is fast an FP64 state
is the simpler choice. Of the less expensive alternatives in Table~\ref{tab:remedies}, only
compensated stages remove the $\Delta t^{-1}$ term, and only for an integrator that is already in
increment form or has been rewritten in it.

The two-float state places the growing term beneath the floor $R_L$ but leaves the floor itself, which is largest where a
small signal depends on the derivative of a large stored value, as for the pressure gradient in
low-amplitude acoustics. In such cases a perturbation state $W'=W-W_{\mathrm{ref}}$ lowers the floor,
because differences such as the pressure gradient are then formed from $p'$ rather than from $p$.
Although reference-state storage is not a substitute for the two-float state
(Section~\ref{sec:refstate}), we recommend it in combination with the two-float state where FP32
accuracy is most critical.

\subsection{Limitations and open questions}\label{sec:implicit}

The analysis does not depend on a particular precision. Two FP64 numbers push the FP64 turnaround
beneath the floor just as two FP32 numbers push the FP32 turnaround beneath it
(Section~\ref{sec:osc}), and since the analysis
is expressed in terms of the unit roundoff of the storage format, an FP16 state with an FP8
right-hand side would follow the same design. The two-float form is useful where the higher
precision is slow or unavailable, as FP64 is on many current accelerators.

The production results were obtained with one solver, using an explicit SSP scheme, for one class
of flows. The wing campaign consists of a single run per build at one flow condition and on one
decomposition, with the two-float run, rather than a repeated FP64 run, serving as the control for
run-to-run variability. The model problems extend the mechanism and the remedy to other physics and
to BDF2. Implicit unsteady methods take larger steps, and the bias is proportional to the number of
steps, so such codes should be less exposed at equal accuracy; their nonlinear residual, however,
cannot be resolved below $\ulp(\bar u)$ with an FP32 state. A pair lowers that limit only if the
time-derivative part of the residual, such as $(3u^{n+1}-4u^{n}+u^{n-1})/(2\Delta t)$ in dual time
stepping, is formed from the pairs, while the spatial operator may still read the leading float. Steady solvers involve no accumulation over steps and are already well
served by mixed-precision iterative correction of a higher-precision state \cite{kelley2022}.

\section{Conclusions}\label{sec:conclusions}

Rounding the stored state is the one rounding error in a time integration that grows as the step
is refined. It acts as a forcing proportional to $1/\Delta t$ whose magnitude is set by the unit in
the last place of the stored mean rather than by the signal, and whose sign is maintained by the
size of the increments or by the rounding of the scheme's coefficients, so that it accumulates as a
bias rather than as a random walk, and refining the step past a turnaround increases the error at
any precision. The ratio of the ulp of the stored mean to the fluctuation of interest indicates in
advance which variables are affected. In wall-modeled large-eddy simulations of a high-lift wing,
the resulting bias reaches $7.7\%$ of the lift on the finest grid and roughly doubles with each
refinement.

Only the stored state requires higher precision. Storing it as two FP32 numbers, updated with
error-free transformations while all right-hand-side evaluations remain in FP32, restores
convergence at the design order down to the right-hand-side floor in every model problem, and brings
the integrated loads of the production simulations to within $0.6\%$ of the FP64 reference on every
grid. Nondimensionalization, reference-state storage and expression-level promotion to FP64 leave the
$1/\Delta t$ term in place, and compensated summation removes it only in increment form, whereas
the two-float state works in the usual form of every integrator tested, from SSPRK(3,3) and
classical RK4 to low-storage SSPRK(10,4) and BDF2. Because the error arises wherever a stored state
is advanced by many small updates, the same construction applies across unsteady scientific
computing, including weather, climate and ocean models, molecular dynamics, astrophysics and
structural simulation, and on any hardware. In the solver studied here it needs about a third less
memory than FP64 and runs up to $2.7$ times faster on a GPU with half-rate FP64, and because it
executes no FP64 operations its advantage grows where FP64 is slower, as on many current GPUs. The
precision that time integration needs belongs in the stored state, not in the arithmetic.

\begin{acknowledgements}
The production simulations use the LAVA Cartesian flow solver, a team code developed at NASA Ames
Research Center whose development the author leads; the author thanks the LAVA team for that
shared foundation and for discussions. Computing resources were provided by the NASA Advanced
Supercomputing (NAS) Division.
\end{acknowledgements}

\section*{Declarations}

\paragraph{Funding} This work was supported by the National Aeronautics and Space Administration
(Entry Systems Modeling Project; Transformational Tools and Technologies Project).

\paragraph{Competing interests} The author has no competing interests to declare that are relevant
to the content of this article.

\paragraph{Data and code availability} The data and code cannot be shared publicly. The production
simulations were performed with the LAVA solver, which is not publicly released, and the
model-problem codes have not been through the NASA software release process. Data supporting the
findings are available from the author on reasonable request, subject to NASA release approval.

\paragraph{Author contributions} Michael F.~Barad is the sole author: conceptualization,
methodology, software, validation, formal analysis, investigation, data curation, writing (original
draft, review and editing) and visualization.

\paragraph{Use of generative AI} During the preparation of this work the author used a large
language model to polish the manuscript and to prototype the model-problem codes. After using this
tool, the author reviewed and verified all code and results, edited the content as needed, and takes
full responsibility for the content of the published article.

\raggedright

\end{document}